\documentclass[12pt,xcolor=svgnames]{amsart}

\usepackage{amsmath, amsthm, amssymb, amsfonts, verbatim}
\numberwithin{equation}{section}

\usepackage{xspace,mathtools,bm}
\usepackage{bigints}

\usepackage[svgnames]{xcolor}

\usepackage[pagebackref=true, colorlinks=true, citecolor=blue]{hyperref}

\allowdisplaybreaks

\makeatletter
\renewcommand{\l@subsection}{\@tocline{2}{0pt}{2em}{3em}{}}
\renewcommand{\l@subsubsection}{\@tocline{3}{0pt}{4em}{4em}{}}
\makeatother

\hypersetup{
    colorlinks=true,
    citecolor=red,
    linkcolor=blue,
    filecolor=magenta,      
    urlcolor=cyan,
    pdfpagemode=FullScreen,
    }
\theoremstyle{plain}
\newtheorem{theorem}{Theorem}[section]

\newtheorem{lemma}{Lemma}[section]
\newtheorem{proposition}{Proposition}[section]
\theoremstyle{definition}
\newtheorem{definition}{Definition}[section]
\newtheorem{remark}{Remark}[section]

\newcommand{\R}{\mathbb{R}}
\newcommand{\eps}{\varepsilon}
\newcommand{\abs}[1]{\mid\!#1\!\mid}
\newcommand{\norm}[2]{\left\lVert#1\right\rVert_{#2}}
\newcommand{\eqdef}{:=}

\newcommand{\Obsolete}[1]{
    }

\DeclareMathOperator*{\esssup}{ess\,sup}
\DeclareMathOperator\supp{supp}

\DeclareFontFamily{U}{mathx}{}
\DeclareFontShape{U}{mathx}{m}{n}{<-> mathx10}{}
\DeclareSymbolFont{mathx}{U}{mathx}{m}{n}
\DeclareMathAccent{\widehat}{0}{mathx}{"70}
\DeclareMathAccent{\widecheck}{0}{mathx}{"71}

\date{\today}
\usepackage[english]{babel}

\title[SQG ill-posedness $W^{k,1}(\R^2)$]{Strong ill-posedness of the Surface Quasi-Geostrophic Equations in $W^{k,1}(\R^2)$}

\author{Elaine Cozzi}
\address{Oregon State University, Department of Mathematics}
\email{elaine.cozzi@oregonstate.edu}

\author{Nicholas Harrison}
\address{Simon Fraser University, Department of Mathematics}
\email{nick\_harrison@sfu.ca}

\author{Zachary Radke}
\address{Oregon State University, Department of Mathematics}
\email{radkeza@oregonstate.edu}
\begin{document}

\begin{abstract}
   We establish ill-posedness of the surface quasi-geostrophic equations in Sobolev spaces $W^{k,1}(\mathbb{R}^2)$ for every integer $k \geq 3$. This includes the scaling critical space $W^{3,1}(\mathbb{R}^2)$. We build on ideas from \cite{BGZ26}; in particular, we use a theorem of Bonami and Poornima to construct a basic building block, or {\em packet}, which obeys a linear evolution approximating that of the perturbation of the background solution to SQG. The Bonami-Poornima theorem ensures that this packet has both small initial $W^{k,1}$ norm and large $W^{k,1}$ norm at positive times. We then implement a process known as {\em many packet dilution}, which allows one to construct an initial perturbation as a sum of sufficiently many packets to dilute the effect of the nonlinear term. This initial perturbation leads to the desired norm inflation.   
   
To establish instantaneous blow-up of the $W^{k,1}(\R^2)$ norm, we also prove a local well-posedness result in the critical Sobolev-Lorentz space $W^{2,(2,1)}(\R^2)$. This result gives a uniform lower bound on the existence time for solutions with bounded $W^{k,1}(\R^2)$ data, even though the solutions' $W^{k,1}(\R^2)$ norms may become arbitrarily large in arbitrarily short time.

\end{abstract}
\maketitle
\tableofcontents

\section{Introduction}
The inviscid surface quasi-geostrophic equations on $\R^2$ are given by
\begin{equation}
\label{SQG} \tag{$SQG$}
\begin{cases}
\partial_t\theta+(u\cdot\nabla)\theta=0 \quad&\text{on } (0,T)\times\R^2,\\
u=\nabla^{\perp}(-\Delta)^{-\frac{1}{2}}\theta = R^\perp \theta&\text{on } [0,T)\times\R^2,\\
\theta\vert_{t=0}=\theta_0 & \text{on } \R^2.
\end{cases}
\end{equation}
The surface quasi-geostrophic equations serve as a model for atmospheric fluid flows and as a two-dimensional analogue to the three-dimensional Euler equations \cite{CMT}.  The scalar quantity $\theta$ denotes the temperature of the fluid, which is transported by a divergence-free velocity field $u$.

Many works in the last three decades have addressed well-posedness for \eqref{SQG} in various function spaces. In this work, we investigate the question of well-posedness of \eqref{SQG} in $L^1$-based Sobolev spaces $W^{k,1}(\R^2)$ for $k\geq 3$. We are motivated by several recent works on well/ill-posedness of fluid equations in $L^1$ based spaces. In \cite{CH}, the first and second authors consider the vorticity formulation of the 2D Euler equations and show short-time persistence of $W^{2,1}(\R^2)$ regularity of the vorticity.  The proof of short-time existence relies on an a priori estimate of the velocity gradient of the form
$$ \| \nabla u \|_{L^{\infty}} \leq C\| \omega \|_{W^{2,1}}.$$
Notably, this a priori control of the velocity gradient fails to hold when vorticity lies in the other scaling critical Sobolev spaces $W^{s,p}(\R^2)$, $sp=2$, $p\in (1,\infty)$. Indeed, in the seminal work \cite{BL}, Bourgain and Li utilize this lack of control of the velocity gradient to construct initial vorticity in $W^{s,p}(\R^2)$, $sp=2$, $p\in (1,\infty)$, which exhibits instantaneous blow-up in the $W^{s,p}(\R^2)$ norm.  

In \cite{DJ}, Danchin and Jeong extend the short-time result of \cite{CH} by demonstrating global-in-time persistence of $W^{2,1}(\R^2)$ regularity of the vorticity. However, due to the unboundedness of the Riesz transforms on $L^1$, the results \cite{CH, DJ} fail to imply persistence of $W^{3,1}(\R^2)$ regularity of the velocity. In fact, in a very recent preprint \cite{BGZ26}, Ban, Guo, and Zhang show that the 2D Euler equations in velocity form are strongly ill-posed in $W^{3,1}(\R^2)$; that is, there exists initial divergence free velocity $u_0\in W^{3,1}(\R^2)$, such that the solution $u$ does not belong to $L^{\infty}([0,T]; W^{3,1}(\R^2))$ for any $T>0$. The main novel ingredient of their proof is the use of the Bonami-Poornima theorem, recalled in Lemma \ref{Bonami Poornima Theorem} below, which characterizes the discontinuous Fourier multipliers on $L^1$-based spaces. Building on the work \cite{BGZ26}, we show that for any $k\geq 3$, \eqref{SQG} is strongly ill-posed in $W^{k,1}(\R^2)$. Notably, for these particular spaces, we have a priori control of the velocity gradient:
$$ \| \nabla u \|_{L^{\infty}} \leq C\|\theta\|_{W^{k,1}}.$$ In particular, in what follows we construct a weak solution to \eqref{SQG} from $\theta_0\in W^{k,1}(\R^2)$ which exhibits instantaneous blow-up in $W^{k,1}$ norm, but for which $\nabla u$ remains bounded for (at least) short time (see Appendix \ref{LWP} for details). For this reason, the methods of Bourgain and Li cannot be used to show norm inflation. Our strategy, which we describe in detail below, more closely follows that of \cite{BGZ26}. 

\subsection{Main Results} We now state our main results.
\begin{theorem}[$W^{k,1}$ Norm Inflation]
\label{norm inflation key}
Let $\overline{\theta}$ be a smooth solution to \eqref{SQG} on $[0,\tau]$ with $\overline{\theta}_0 \in C_c^{\infty}(\R^2)$ non-constant. Then, for every $\varepsilon,M>0$, and integer $k \geq 3$, there exists $\theta_0 \in C_c^{\infty}(\R^2)$ whose unique solution to \eqref{SQG} exists through $\tau$,
\begin{equation}
\label{norm inflation key equation}
\norm{\theta_0-\overline{\theta}_0}{W^{k,1}}<\varepsilon, \quad \text{and}\quad \norm{\theta(\tau)-\overline{\theta}(\tau)}{W^{k,1}} >M.
\end{equation}
Moreover, if $L$ is such that $\supp \overline{\theta_0} \subseteq B_L(0)$, then $\theta_0$ can be chosen to satisfy $\supp \theta_0 \subseteq B_L(0)$. 
\end{theorem}
\begin{theorem}[Instantaneous $W^{k,1}$ Blow Up] 
\label{BLOW UP}Let $k \geq 3$ be an integer, $\overline{\theta}_0 \in C_c^{\infty}(\R^2)$ be non-constant, and $\eps>0$. Then, there exists a compactly supported $\theta_0 \in W^{k,1}(\R^2)$ and $T>0$ such that
\begin{equation}
\label{blow up close data}
\norm{\theta_0-\overline{\theta}_0}{W^{k,1}} <\eps,
\end{equation}
and the corresponding unique weak solution $\theta$ to \eqref{SQG} on $[0,T]$, constructed in Appendix \ref{LWP}, satisfies
\begin{equation}
\label{instant blow up}
\esssup_{0<t<T'} \norm{\theta(t)}{W^{k,1}}=\infty, \quad \text{for every } 0<T'\leq T.
\end{equation}
\end{theorem}

\begin{remark} 
While we expect ill-posedness of \eqref{SQG} to hold in $W^{k,1}(\R^2)$ for $k=1,2$, our argument does not directly extend to these values of $k$. To illustrate the obstruction, consider a localized initial perturbation of the form 
$$
\eta_0(x) = ar^{k-2}f\left(\frac{x-q}{r}\right),
$$
where $f \in C_c^{\infty}(\R^2)$, $q \in \R^2$, $a>0$, and $0<r\leq 1$. From scaling
$$
\|\eta_0\|_{\dot{W}^{k,1}} = a\|f\|_{\dot{W}^{k,1}}, \quad \|\nabla R^{\perp}\eta_0\|_{L^{\infty}} = ar^{k-3}\|\nabla R^{\perp}f\|_{L^{\infty}}.
$$ 
For $ k \geq 3$, the velocity gradient is therefore bounded by a constant times $a$, uniformly as $r \to 0$. When $k=1,2$, however, concentration increases the velocity gradient. This prevents the same estimates later used from providing uniform control of the perturbations' self-interaction as $r \to 0$.

Special structure can sometimes provide cancellations that are not captured by this scaling estimate. For example, a radial profile $f$ satisfies $R^{\perp}f\cdot \nabla f=0$. However, the Bonami-Poornima theorem supplies us with the profile to obtain norm inflation, but does not provide us with the symmetry needed to exploit this cancellation. Furthermore, even a radially symmetric initial profile need not remain radially symmetric under the linear evolution considered in Section \ref{norm inflation section}. Hence, extending our argument to $k=1,2$ would require additional estimates or a different choice of perturbation.
\end{remark}

\subsection{Literature Review}
The study of fluid equations in Sobolev spaces has a long and deep history; we do not attempt to give a complete overview, but discuss some landmark results to contextualize the significance of the present work. We begin with the Euler equations in vorticity form. For $1<p<\infty$, following an argument from \cite{KP}, one obtains an a priori estimate of the form 
\begin{equation}
\label{2d vorticity apriori}
\frac{d}{dt}\norm{\omega(t)}{W^{s,p}} \leq C\norm{\nabla u}{L^{\infty}}\norm{\omega(t)}{W^{s,p}},
\end{equation}
where $u=\nabla^{\perp}(-\Delta)^{-1}\omega$. In order to close \eqref{2d vorticity apriori}, for $s>\frac{d}{p}$ and $1<p<\infty$, one can apply the Sobolev embedding theorem and the boundedness of Riesz transforms on $L^p(\R^d)$  to write $\|\nabla u\|_{L^{\infty}} \leq C\|\omega\|_{W^{s,p}}$, yielding short-time existence after a suitable approximation argument. The work of Beale, Kato, and Majda in \cite{BKM} then provides a criterion for the continuation of the solution, namely $$\limsup_{t \to T} \|\omega(t)\|_{W^{s,p}} = \infty \iff \int_0^T \|\omega(t)\|_{L^{\infty}}\,dt = \infty.$$
Consequently, in the two-dimensional case, the conservation of the $L^\infty(\R^2)$ norm of the vorticity yields global well-posedness in $W^{s,p}(\R^2)$, provided that $sp>2$ and $1<p<\infty$. On the other hand, when $sp \leq d$ and $p>1$, the inequality $\norm{\nabla u}{L^{\infty}} \leq C\norm{\omega}{W^{s,p}}$ fails to hold. For this reason, we refer to $W^{s,p}(\R^d)$ as supercritical, critical, or subcritical when $sp>d$, $sp=d$, and $sp<d$, respectively. 

The analysis of the Euler equations in Sobolev spaces was subsequently extended to the more general Triebel-Lizorkin and Besov spaces, see, for example, \cite{ChaeBesov,ChaeTriebel, Vishik98, PakPark04,PakPark13,HwangPak,HarrisonRadke}. As with the Sobolev spaces, in supercritical Triebel-Lizorkin and Besov spaces (i.e. those for which a priori control of the velocity gradient is available), standard techniques establish local well-posedness when $d \geq 2$ and global well-posedness when $d=2$. For the critical spaces, in the two-dimensional setting, Vishik \cite{Vishik98} establishes global well-posedness in $B^0_{\infty,1}(\R^2)$ using the a priori estimate $\| \nabla u \|_{L^{\infty}}\leq C\|\omega\|_{B^0_{\infty,1}}$, para-differential calculus, and the transport structure of the vorticity equation. Analogously, Chae \cite{chae2003euler} obtains global well-posedness with vorticity in the critical Triebel-Lizorkin spaces $F^2_{1,q}(\R^2)$. 

Without a priori control of the velocity gradient, the situation is considerably more delicate. In the seminal work \cite{BL}, Bourgain and Li utilize this lack of a priori control to establish strong ill-posedness of the two-dimensional Euler equations (in vorticity form) in critical Sobolev spaces $W^{s,p}(\R^2)$ for $sp = 2$ and $1<p<\infty$. To obtain ill-posedness, Bourgain and Li construct a small initial vorticity whose associated flow develops large Lagrangian deformation, which can be achieved using a key lemma from \cite{KS}. They then introduce a highly oscillatory perturbation and show that stretching of the particle-trajectory map can produce norm inflation.

As with the Euler equations, the question of well-posedness for \eqref{SQG} in various function spaces has received considerable attention over the past three decades. This line of investigation began with the work of Constantin-Majda-Tabak \cite{CMT}, who establish local existence and uniqueness for initial data in $H^k(\R^2)$, with integer $k \geq 3$, together with a continuation criterion involving the time integral of $\|\nabla \theta(t)\|_{L^{\infty}}$. Chae \cite{Chae03} subsequently develops local well-posedness and blow-up criteria in Triebel-Lizorkin spaces. In \cite{WangJia}, Wang and Jia establish local well-posedness in the critical Besov spaces $B^{2/p+1}_{p,1}(\R^2)$, for $1 \leq p \leq \infty$. In particular, their results include $B^3_{1,1}(\R^2)$, but do not include the larger space $W^{3,1}(\R^2)$.

Several recent results have demonstrated ill-posedness of \eqref{SQG} when a priori control of $\|\nabla\theta\|_{L^{\infty}}$ is unavailable.
Jeong and Kim \cite{JeongKim} use this lack of control to establish strong ill-posedness in $H^2(\R^2)$ for \eqref{SQG} by a contradiction argument. They assume that such a solution belongs to $L^{\infty}([0,T];H^2(\R^2))$, and then develop an analogue of the Kiselev-\v{S}ver\'{a}k estimate, which states that if the initial scalar $\theta$ is odd-odd then
\begin{equation}
\label{SQG velocity approx}
\frac{u_1}{x_1} \approx \int_{[2x_1,1]\times [0,1]} \frac{y_1y_2}{|y|^5}\theta(y)\,dy, \quad |x| \ll 1,
\end{equation}
providing an approximation of the velocity gradient at the origin. They use \eqref{SQG velocity approx}, a hypothetical $H^2$-bound of $\theta$, and careful tracking of dyadic packets to show that such a solution cannot exist for any $t>0$.

Independently, in \cite{CordobaZoroa}, C\'{o}rdoba and Mart\'inez-Zoroa also consider the problem of ill-posedness for \eqref{SQG} in $H^2(\R^2)$, but take a different approach than Jeong and Kim. They write $x \in \R^2$ in polar coordinates as $x=(r\cos(\alpha),r\sin(\alpha))$ and consider initial data of the form 
$$
\theta_0^J(x)=\sum_{j=1}^J\theta_{0,j}^J(x), \quad \theta_{0,j}^J= c\frac{f(b^{-j}r)b^j\sin(2\alpha)}{j},
$$
where $f \in C_c^{\infty}((0,\infty))$ is supported in some annulus. Note that for each $J$, $\theta_0^J \in C_c^{\infty}(\R^2)$ and hence generates a unique smooth solution. Furthermore, the initial data $\theta_0^J$ satisfies 
$$
\|\theta_0^J\|_{\dot{H^2}} \lesssim c \quad\text{and}\quad \|\nabla u_0(0)\|_{L^{\infty}} \sim c\log(J).
$$
The authors analyze the evolution of each $\theta_{0,j}^J$ and send $b \to 0^+$ (which spatially separates each $\theta_{0,j}^J$). However, increasing the separation between scales does not minimize pairwise interaction. To work around this, they show that the larger annuli generate a leading hyperbolic strain on each smaller annulus, which they retain explicitly while controlling the remaining interaction errors to ensure that the solution persists long enough for substantial deformation to occur, opening the door to introduce a highly oscillatory perturbation, yielding norm-inflation.

The well/ill-posedness results described above for the Euler equations and \eqref{SQG} address $L^p$-based spaces when $p\in (1,\infty)$; however, the $L^1$-based spaces exhibit strikingly different behavior. For example, unlike in the critical spaces $W^{s,p}(\R^2)$, $sp=2$, the two-dimensional Euler equations (in vorticity form) are globally well-posed in $W^{2,1}(\R^2)$. Indeed, using Lorentz spaces to utilize finer embeddings, the authors in \cite{CH} note that $\|\nabla u\|_{L^{\infty}} \leq C\|\omega\|_{W^{2,1}}$, establishing the necessary a priori estimate of the velocity gradient to obtain local well-posedness. Global well-posedness with vorticity in $W^{2,1}(\R^2)$ is subsequently established by Danchin and Jeong in \cite{DJ}. 

In analogy with the Euler equations, the velocity and scalar in \eqref{SQG} satisfy the estimate $\|\nabla u\|_{L^{\infty}} \leq C\|\theta\|_{W^{3,1}}$, which naturally leads to the question of local well-posedness of \eqref{SQG} in $W^{3,1}(\R^2)$. With this a priori control of the velocity gradient, one may expect local well-posedness to hold. In particular, the methods of Bourgain and Li cannot establish ill-posedness, as these rely on large velocity gradients with small norm of the scalar. At the same time, the endpoint estimate above falls short of controlling the highest-order non-local term $\|D^3u \cdot \nabla \theta\|_{L^1}$. To establish ill-posedness, a novel approach appears necessary. 

In a very recent preprint \cite{BGZ26}, the authors establish strong ill-posedness of the two-dimensional Euler equations in $W^{3,1}(\R^2)$. We make two remarks on their work. First, in their setting, $\|\nabla u\|_{L^{\infty}} \leq C\|u\|_{W^{3,1}}$, which presents the same obstacle to showing ill-posedness that we face in this paper. Second, the $W^{2,1}(\R^2)$ norm of the vorticity is bounded globally in time; hence, the mechanism for ill-posedness must be a consequence of how the velocity is reconstructed from the vorticity. Their proof begins by placing a small, localized perturbation in a smooth background flow which stretches particles in one direction and compresses them in another. Near a background trajectory $X(t)$, they approximate the perturbation's vorticity by the linear equation
$$
\partial_t \eta + (\nabla \overline{u}(t,X(t))y)\cdot \nabla \eta=0, \quad \eta(0,y)=-\Delta f(y),
$$
where $f$ is a smooth, compactly supported initial stream function corresponding to the perturbation. After solving this equation, they recover the approximate perturbation velocity by $v=\nabla^{\perp}(-\Delta)^{-1}\eta$. In particular, if $B(t)$ is the solution to
$$
B'(t) = \nabla \overline{u}(t,X(t))B(t), \quad B(0)=\text{Id}_{2 \times 2},
$$
then 
$$
\hat{v}(t,\xi)=2\pi i\frac{\xi^{\perp}}{|\xi|^2}|B(t)^T\xi|^2\hat{f}(B(t)^T\xi),
$$
so by a change of variables
$$
v(t,B(t)y)=B(t)\nabla^{\perp}[m_t(D)f](y), \quad m_t(\xi)= \frac{|\xi|^2}{|B(t)^{-T}\xi|^2}.
$$
One easily checks that $m_t(D)$ is a smooth Fourier multiplier which is homogeneous of degree $0$, and under more technical assumptions (namely that the background flow is anisotropic), is also non-constant, allowing the authors to apply the Bonami–Poornima theorem. In particular, for any $L>0$, there exists a stream function $f \in C_c^{\infty}(\R^2)$ which is mean zero, $\|f\|_{\dot{W}^{4,1}} \leq 1$, and $\|m_t(D)f\|_{\dot{W}^{4,1}}>L$. The remaining work is to show that the norm inflation of the approximate perturbation velocity is not destroyed by the nonlinear dynamics.

In this work, we apply a similar strategy to that of \cite{BGZ26} to establish strong ill-posedness of \eqref{SQG} in $W^{k,1}(\R^2)$ for $k \geq 3$, namely, the use of the Bonami-Poornima theorem to construct initial data. However, our approach differs from that of \cite{BGZ26} in several ways. First, for the Euler equations, the mechanism responsible for norm inflation arises from the anisotropic deformation of the vorticity and subsequent recovery of the velocity through the Biot-Savart law. In contrast, linearizing \eqref{SQG} around a smooth background solution produces an additional nonlocal term $\nabla \overline{\theta}\cdot R^{\perp}\eta$, which gives rise to a non-constant (provided that $\nabla \overline{\theta}_0$ is non-zero), smooth, and homogeneous of degree $0$ symbol, even in the absence of anisotropic deformation by the background flow. Hence, our instability mechanism is driven directly by the nonlocal interaction of the perturbation and the background scalar, rather than recovering the velocity from the vorticity. Moreover, in \cite{BGZ26} one has $\nabla u \in L^{\infty}_{loc}([0,\infty);L^{\infty}(\R^2))$ from $\omega \in L^{\infty}_{loc}((0,\infty);W^{2,1}(\R^2))$. We do not have an analogous result; thus, when establishing instantaneous blow-up, we must also prove local well-posedness for \eqref{SQG} in the Sobolev-Lorentz space $W^{2,(2,1)}(\R^2)$ (see Appendix \ref{LWP}). In particular, the embedding $W^{3,1}(\R^2) \hookrightarrow W^{2,(2,1)}(\R^2)$ allows us to obtain a common time of existence for a sequence of smooth solutions with uniformly bounded initial $W^{3,1}(\R^2)$ norms. We describe our proof strategy in greater detail below.

\subsection{Proof Strategy} We first describe the proof of norm inflation. Let us write $\eta = \theta-\overline{\theta}$ as a perturbation of $\overline{\theta}$, where $\theta,\overline{\theta}$ are as in Theorem \ref{norm inflation key}. Then subtracting the equations satisfied by $\theta$ and $\overline{\theta},$ we obtain
$$\partial_t\eta+\overline{u}\cdot\nabla \eta +R^\perp\eta\cdot\nabla \overline{\theta}+(R^\perp\eta)\cdot\nabla \eta =0.$$

Our goal is to show $\eta$ can be constructed so that $\eta|_{t=0}$ is small in $W^{k,1}(\mathbb{R}^2)$ but $\|\eta(t)\|_{W^{k,1}}$ can be made very large at small times $t>0$. Each of the terms above plays a separate role: $(\partial_t+\overline{u}\cdot \nabla) \eta $ transports the perturbation by the background velocity $\overline{u}$; $R^\perp \eta\cdot \nabla \overline{\theta}$ contains the velocity created by the perturbation which then transports the background scalar $\overline{\theta}$; and $R^\perp \eta\cdot\nabla \eta$ is the nonlinear self-advection of the scalar, which we must keep small throughout to identify the growth of $\eta$ with the linear evolution of the other terms.

In the first step, we replace the background velocity in the above equation by a linear Taylor approximation and the gradient of the background scalar by its value on a given trajectory. If we now rewrite the equation in the coordinates given by the deformation of the scalar, neglecting the nonlinear term, we obtain a linear evolution corresponding to multiplication by a Fourier multiplier which is non-constant and homogeneous of degree $0$. We then apply the theorem of Bonami and Poornima, recalled in Lemma \ref{Bonami Poornima Theorem} below, to obtain growth of the solution to the linear equation.

The validity of the construction then hinges on showing that the construction itself can be made to minimize the effect of the nonlinear term $R^\perp\eta\cdot\nabla \eta $ while also maintaining the amplification of the $W^{k,1}(\mathbb{R}^2)$ norm, which is done through a process referred to as \textit{many-packet dilution} in \cite{BGZ26}. Specifically, with $\varepsilon>0$ being the total size of the initial perturbation, we construct a sum with $J$ copies of the packet obtained from the Bonami-Poornima theorem, each of size $a=\varepsilon/J$. We show that this choice of scale $a$ allows us to estimate the error due to the nonlinear term by an upper bound which can be made arbitrarily small by choosing $J$ large. The remaining error is controlled by rescaling the spatial concentration of the packet. This idea is further explained in Remark \ref{Remark dilution}, which references the specific contributions made by the different approximation choices to the main $W^{k,1}(\mathbb{R}^2)$ lower bound \eqref{LB1}.

For the proof of instantaneous blow-up, given non-constant $\overline{\theta}_0 \in C_c^{\infty}$, we create a sequence of initial data $\theta_0^{(n)}$ by inductively applying the norm-inflation result. To do this, given $\theta_0^{(n-1)}$, we create $\theta_0^{(n)}$ and a sequence of times $t_n \to 0$ so that 
$$
\norm{\theta^{(n)}(t_n)}{W^{k,1}} > 4n \quad \text{and}\quad \sum_{j=1}^{\infty} \norm{\theta_0^{(j)}-\theta_0^{(j-1)}}{W^{k,1}}<\infty.
$$
Consequently, the $\theta_0^{(n)}$ reside in a ball of uniform radius around $\overline{\theta}_0$ in $W^{2,(2,1)}$. Applying the result from Appendix \ref{LWP}, we obtain a sequence of unique weak solutions $\theta^{(n)} \in L^{\infty}([0,T];W^{2,(2,1)}(\R^2))$ which are uniformly bounded. In particular, we see that $\nabla \theta^{(n)} \in L^1([0,T];L^{\infty}(\R^2))$ is uniformly bounded, allowing us to establish ``stability" in $L^2(\R^2)$.

This stability estimate lets us preserve the large norms created at each stage under sufficiently small further changes to the initial data. To see this, we first localize and mollify $\theta^{(n)}(t_n)$, and then analyze its difference with another given localized and mollified solution at time $t_n$. Pushing all derivatives onto the mollifier and using the Cauchy-Schwarz inequality lets us control the difference between these solutions by their $L^2(\R^2)$ difference, which we can bound by the $L^2(\R^2)$ difference of the initial data via stability. This control allows us to choose $\rho_n>0$ such that any initial data within $\rho_n$ of $\theta_0^{(n)}$ in $W^{k,1}(\R^2)$ produces a solution which, when localized and mollified, has a $W^{k,1}(\R^2)$ norm at time $t_n$ of at least $2n$.

We then make all subsequent perturbations sufficiently small so that their total size is less than $\rho_n$. We show that the sequence of initial data is Cauchy in $W^{k,1}(\R^2)$, hence it converges to some initial data $\theta_0 \in W^{k,1}(\R^2)$. If all of the perturbations have support contained in a ball $B_L(0)$ which contains $\supp \overline{\theta}_0$, then $\supp \theta_0$ will also be contained in $\overline{B_L(0)}$. Then, from the continuity of the $W^{k,1}(\R^2)$ norm, $\theta_0$ lies within $\rho_n$ of each $\theta_0^{(n)}$, showing that its localized and mollified solution inherits a lower bound of $2n$ on its $W^{k,1}(\R^2)$ norm at time $t_n$. Since each localized mollified functional is bounded above by the full $W^{k,1}(\R^2)$ norm, the lower bounds for these functionals imply the required norm growth for instantaneous blow-up.

\subsection{Notational Conventions}
In what follows, $B_r(x)$ denotes the ball of radius $r>0$ about the point $x \in \R^d$. For ease of notation, we set $B_r(0)=B_r$. We let $|A|$ denote the Lebesgue measure of a set $A \subseteq \R^d$. Below, $C>0$ denotes an arbitrary constant which may change from line to line. If the constant depends on a set of parameters $\{\lambda_1,\cdots,\lambda_n\}$ and we aim to stress this dependence, we write $C_{\lambda_1,\cdots,\lambda_n}$. If there exist constants $C_1$ and $C_2$ such that $A \leq C_1B$ and $B \leq C_2A$, we write $A \sim B$. Given two functions $f,g:(0,\infty) \to \R$, with $g(r) \neq 0$ for $r$ sufficiently small, we write 
$$
f(r)=o(g(r)),
$$
if $\lim_{r \to 0^+} \frac{f(r)}{g(r)}=0$. Similarly, we write $$f(r)=\mathcal{O}(g(r)),$$ if there exist constants $C>0$ and $r_0>0$ such that $|f(r)| \leq C|g(r)|$ for all $r<r_0$. For $f \in \mathcal{S}(\R^d)$, we  use the convention
$$\mathcal{F}(f)(\xi)=\hat{f}(\xi) := \int_{\mathbb{R}^d} f(x) e^{-2\pi i x \cdot \xi}\, dx $$ to denote the Fourier transform of $f$. For $x \in \R^2$ we define
$$
x^{\perp}\eqdef (-x_2,x_1), \quad \nabla^{\perp}\eqdef (-\partial_2,\partial_1), \quad R^{\perp}\eqdef \nabla^{\perp}(-\Delta)^{-\frac{1}{2}}.
$$
With our convention of the Fourier transform, note that for smooth $f$
$$
\widehat{R^{\perp}f}(\xi) = i\frac{\xi^{\perp}}{|\xi|}\hat{f}(\xi), \quad \xi \neq 0.
$$
Additionally, given $f\in \mathcal{S}(\R^d)$ and $m \in C^{\infty}(\R^d \setminus \{0\})$, we define 
$$m(D)f \eqdef \int_{\R^d}m(\xi) \hat{f}(\xi)e^{2\pi i x \cdot \xi }\,d\xi.$$
\subsection{Organization of the Paper} In Section \ref{norm inflation section} we establish norm inflation, proving Theorem \ref{norm inflation key}. In Section \ref{instant blow up section} we demonstrate instantaneous blow-up of the $W^{k,1}(\R^2)$-norm, proving Theorem \ref{BLOW UP}. Appendix \ref{appendix energy} contains the energy estimates utilized in Section \ref{norm inflation section}. Finally, in Appendix \ref{LWP}, we show local well-posedness in the Sobolev-Lorentz space $W^{2,(2,1)}(\R^2)$.
\section{Proof of Theorem \ref{norm inflation key}: Norm Inflation}
\label{norm inflation section}
Throughout this section, we let $\bar{\theta}$ and $\tau$ be as in Theorem \ref{norm inflation key}.
\subsection{The Building Block}
In this subsection, we construct the basic building block used in the iterative construction of $\theta_0$ in Theorem \ref{norm inflation key}.
Again, let $\eta =\theta-\overline{\theta}$ denote the perturbation of $\overline{\theta}$. Then, the perturbation $\eta$ satisfies
\begin{equation}
\label{initial perturbation pde}
\partial_t\eta+\overline{u}\cdot\nabla \eta +R^{\perp}\eta \cdot \nabla \overline{\theta}+R^{\perp}\eta \cdot \nabla\eta=0.
\end{equation}Let $q\in \mathbb{R}^2$ and let $X_q(t)$ denote the trajectory of the point $q$ under the background flow $\overline{u}$. That is,
\begin{equation}\label{qtrajectory}\frac{d}{dt}X_q(t) = \overline{u}(t,X_q(t)),\qquad X_q(0)=q.
\end{equation}
To motivate the equation approximating (\ref{initial perturbation pde}), we fix $q \in \R^2$ and assume we have an initial perturbation of the form 
$$
\eta(0,x)=\eta_0(x-q).
$$
We then express the perturbation in moving coordinates. Set
$$
y=x-X_q(t), \qquad \eta_q(t,y)=\eta(t,X_q(t)+y).
$$
Since $X_q(0)=q$, $$\eta_q(0,y)=\eta(0,q+y)=\eta_0(y).$$ 
By the chain rule and \eqref{qtrajectory}
\begin{equation}
\label{chain rule perturbation}
\partial_t\eta(t,x)=\partial_t\eta_q(t,y)-X_q'(t)\cdot \nabla_y\eta_q(t,y).
\end{equation}
If $\mathcal{K}$ is the kernel corresponding to $R^{\perp},$
\begin{equation}
\begin{split}
\label{translated riesz transform}
R^{\perp}\eta(t,x)&= \int_{\R^2}\mathcal{K}(x-z)\eta_q(t,z-X_q(t))\,dz = \int_{\R^2}\mathcal{K}(x-z-X_q(t))\eta_q(t,z)\,dz\\
&= \int_{\R^2}\mathcal{K}(y-z)\eta_q(t,z)\,dz = R^{\perp}\eta_q(t,y).
\end{split}
\end{equation}
Then, substituting \eqref{chain rule perturbation} and \eqref{translated riesz transform} into \eqref{initial perturbation pde}, one finds that the exact translated perturbation $\eta_q$ satisfies
$$
\partial_t\eta_q+(\overline{u}(t,X_q(t)+y)-X_q'(t))\cdot\nabla\eta_q + \nabla \overline{\theta}(t,X_q(t)+y)\cdot R^{\perp}\eta_q + R^{\perp}\eta_q\cdot\nabla \eta_q=0.
$$
Using a Taylor approximation
\begin{align*}
     \overline{u}(t,X_q(t)+y)-X_q'(t) &= \overline{u}(t,X_q(t)+y)-\overline{u}(t,X_q(t)) \\&=\nabla\overline{u}(t,X_q(t))y+\mathcal{O}(|y|^2),  
 \end{align*}
 and 
 \begin{equation*}
     \nabla\overline{\theta}(t,X_q(t)+y) = \nabla\overline{\theta}(t,X_q(t))+\mathcal{O}(|y|).
 \end{equation*}
Motivated by these expansions, we define the building block $\psi_q$ by retaining the leading background terms and omitting the nonlinear self-advection term to obtain
 \begin{equation}
 \label{Packet Equation}
 \begin{cases}
\partial_t\psi_q + (\nabla \overline{u}(t,X_q(t))y)\cdot\nabla\psi_q + R^{\perp}\psi_q \cdot \nabla \overline{\theta}(t,X_q(t))=0,\\
\psi_q\vert_{t=0}=\eta_0.
 \end{cases}
 \end{equation}
 The amplitude and concentration parameters, which we introduce later in Section \ref{norm inflation section}, will allow us to control the errors arising from the above approximations.

The function $\psi_q$ describes the approximate perturbation evolution in coordinates moving along $X_q(t)$. Translating back to physical coordinates gives $\psi_q(t,x-X_q(t))$, whose initial value is $\eta_0(x-q)$. Thus, this change of coordinates recovers the location of the original packet. We now consider the linear model for varying centers $q$, using the same initial profile $\eta_0$ for each center.  To explain why a single profile is sufficient, in Lemma \ref{multiplier lemma} we show that if for some $q_0\in\R^2$, $\nabla \overline{\theta}_0(q_0) \neq 0$, then the solution $\psi_{q_0}$ of \eqref{Packet Equation} corresponds to a non-constant, smooth (away from the origin) Fourier multiplier which is homogeneous of degree zero. In Theorem \ref{Packet Family Theorem}, we then leverage the Bonami-Poornima theorem to construct an initial perturbation that yields large amplification of the norm $\|\psi_{q_0}(\tau)\|_{W^{k,1}}$  at time $\tau$. We proceed to show that if $q$ is sufficiently close to $q_0$, $\|\psi_{q}(\tau)\|_{W^{k,1}}$ can be made large as well.

The following lemma establishes a formula for the solution $\psi_q$ to \eqref{Packet Equation} and prepares us for the use of the Bonami-Poornima theorem.
\begin{lemma}
    \label{multiplier lemma}
    Let 
    \begin{equation}\label{definition of Bq}
    \frac{d}{dt}B_q(t) = \nabla \overline{u}(t,X_q(t))B_q(t),\qquad B_q(0)  = Id_{2\times 2}.
    \end{equation}
    Then $\det B_q = 1$ and the solution $\psi_q$ of (\ref{Packet Equation}) satisfies  
    \begin{equation}\label{ML:psiq solution formula}\psi_q(t,y) =\big[m_q(t,D)\eta_0](B_q(t)^{-1}y),
    \end{equation} where, for $\xi\not=0$,
    \begin{equation}
        \label{ML:multiplier formula}
        m_q(t,\xi) = \exp\left\{-i\big(\nabla\overline{\theta}_0(q)\cdot\xi^\perp\big)\int_0^t\frac{ds}{|B_q(s)^{-T}\xi|}\right\}.
    \end{equation}
    If $\nabla \overline{\theta}_0(q)\not = 0$, then for every $0<t\leq \tau$ this multiplier is smooth away from zero, non-constant, and homogeneous of degree zero.
\end{lemma}
\begin{proof}
Since $\overline{u}$ is divergence free, $
\det B_q(t)=  1,
$
see for example Chapter $2$ in \cite{MB}.
Furthermore, since 
$$
(\partial_t + \overline{u}\cdot \nabla)\nabla \overline{\theta} = (-\nabla \overline{u})^T(\nabla \overline{\theta}),
$$
we obtain 
\begin{equation}
\label{ML: nabla background 1}
\frac{d}{dt} \nabla\overline{\theta}(t,X_q(t))=-\nabla\overline{u}^T(t,X_q(t))\nabla\overline{\theta}(t,X_q(t)).
\end{equation}
Since  $B_q^{-1}B_q=Id_{2\times 2}$, by the product rule we have $(B_q^{-1})'B_q=-B_q^{-1}B_q'$ and hence
$$
(B_q^{-1})' = -B_{q}^{-1}B_q'B_{q}^{-1}=-B_{q}^{-1}\nabla\overline{u}(t,X_q(t))B_q{B_{q}^{-1}} = -B_{q}^{-1}\nabla\overline{u}(t,X_q(t)).
$$
Taking the transpose and multiplying by $\nabla\overline{\theta}_0(q)$ we see that 
\begin{equation}
\label{ML: nabla background 2}
\begin{cases}
\frac{d}{dt}(B_q^{-T}\nabla\overline{\theta}_0(q))=-\nabla \overline{u}^T(t,X_q(t))(B_q^{-T}\nabla\overline{\theta}_0(q)),\\
B_{q}^{-T}\nabla\overline{\theta}_0(q)\vert_{t=0}=\nabla\overline{\theta}_0(q).
\end{cases}
\end{equation}
In particular, comparing \eqref{ML: nabla background 1} and \eqref{ML: nabla background 2}, one finds by uniqueness that 
\begin{equation}
\label{ML: nabla background}
\nabla\overline{\theta}(t,X_q(t))=B_q(t)^{-T}\nabla\overline{\theta}_0(q).
\end{equation}
Evaluating $\psi_q$ at $y=B_q(t)z$, it follows from the multi-variable chain rule and \eqref{Packet Equation} that 
\begin{equation}
\label{ML:preFourier}
\partial_t(\psi_q(t,B_q(t)z)) = -\nabla\overline{\theta}(t,X_q(t))\cdot (R^{\perp}\psi_q)(t,B_q(t)z),\qquad \psi_q(0,B_q(0)z)=\eta_0(z).
\end{equation}
 Hence, setting 
 $$f(t,\xi) = \mathcal{F}_z[\psi_q(t,B_q(t)z)](\xi)$$
 and applying the Fourier transform (in $z$) to \eqref{ML:preFourier}, we have 
\begin{equation*}
\begin{split}
&\partial_tf(t,\xi) = -\nabla\overline{\theta}(t,X_q(t))\cdot \mathcal{F}_z[(R^{\perp}\psi_q)(t,B_q(t)z)](\xi)\\
&\qquad = -\nabla\overline{\theta}(t,X_q(t))\cdot \widehat{R^{\perp}\psi_q}(B_q(t)^{-T}\xi)\\
&\qquad = -i\frac{\nabla \overline{\theta}(t,X_q(t))\cdot\big(B_q(t)^{-T}\xi\big)^\perp}{|B_q(t)^{-T}(\xi)|}f(t,\xi),
\end{split}
\end{equation*}
where we used a change of variables and that $\det B_q(t) = 1$ to obtain the second and third equalities. We conclude that
$$f(t,\xi) = \widehat{\eta_0}(\xi)\exp\left\{-i\bigintssss_0^t\frac{\nabla\overline{\theta}(s,X_q(s))\cdot\Big(B_q(s)^{-T}\xi\Big)^\perp }{|B_q(s)^{-T}\xi|}\,ds\right\}.$$
 Substituting $\eqref{ML: nabla background}$ in the numerator of the integrand, applying the two-dimensional identity $$Bv\cdot (Bw)^\perp = \det(B)v\cdot w^\perp$$
 to the resulting expression,
 and utilizing the fact that $\det(B_q(s)^{-T})=1$, we obtain 
\begin{equation}
\label{ML:postFouriersolution}
f(t,\xi) =  \widehat{\eta_0}(\xi)\exp\left\{-i(\nabla\overline{\theta}_0(q)\cdot \xi^{\perp})\int_0^t \frac{ds}{|B_q(s)^{-T}\xi|}\right\}.
\end{equation}
Inverting the Fourier transform now, we have
\begin{equation}
\label{ML:postinverse}
\psi_q(t,B_q(t)z) = \big[m_q(t,D)\eta_0\big](z)
\end{equation}
where $m_q$ is as in $\eqref{ML:multiplier formula}$. Setting $y = B_q(t)z$ yields (\ref{ML:psiq solution formula}).

To complete the proof of Lemma \ref{multiplier lemma}, it remains to show that when $\nabla\bar{\theta}_0(q) \neq 0$, $m_q$ satisfies the desired properties. Assume $\nabla\overline{\theta}_0(q) \neq 0$ and fix $t\in (0,\tau]$. We first show that $m_q$ is non-constant. We define $$\xi_1 = \frac{\nabla \overline{\theta}_0(q)}{|\nabla \overline{\theta}_0(q)|}.$$
It immediately follows that $m_{q}(t,\xi_1) =1$. Next, we rotate $\xi_1$ by an angle $\alpha$. Hence, we define
$$
\xi_{\alpha} = \cos(\alpha)\xi_1 + \sin(\alpha)\xi_1^{\perp}.$$
We claim that for $\alpha>0$ sufficiently small, $m_q(t,\xi_{\alpha})\neq 1$.  To see this, first note that  $\xi_{\alpha}^{\perp} = \cos(\alpha)\xi_1^{\perp} - \sin(\alpha)\xi_1$ so that 
\begin{equation}\label{xialphaconstruction}
\nabla \overline{\theta}_0(q) \cdot \xi_{\alpha}^{\perp}=-|\nabla \overline{\theta}_0(q)|\sin(\alpha).
\end{equation}
We set $$K=\max\{1,\sup_{0 \leq s \leq t} \|B_{q}(s)^T\|\},$$ which we know is finite due to smoothness of $\overline{\theta}$. Then $1=|B_{q}(s)^TB_{q}(s)^{-T}\xi_{\alpha}| \leq K|B_{q}(s)^{-T}\xi_{\alpha}|.$ It follows that 
\begin{equation}\label{integralKbound}
0<\int_0^t \frac{ds}{|B_{q}(s)^{-T}\xi_{\alpha}|} \leq Kt.
\end{equation}
We then choose $\alpha$ satisfying
\begin{equation}\label{alphabound}
0<\alpha<\min \left\{\frac{\pi}{2}, \frac{\pi}{Kt |\nabla\overline{\theta}_0(q)|}\right\}.
\end{equation}
Since $0<\sin\alpha \leq \alpha$,  (\ref{xialphaconstruction}), (\ref{integralKbound}), and (\ref{alphabound}) yield 
$$
0 < \left|\nabla \overline{\theta}_0(q)\cdot \xi_{\alpha}^{\perp} \int_0^t \frac{ds}{|B_{q}(s)^{-T}\xi_{\alpha}|}\right| \leq |\nabla \overline{\theta}_0(q)|Kt\alpha < \pi.
$$
Consequently,
$
m_{q}(t,\xi_{\alpha}) \neq 1 = m_{q}(t,\xi_1)
$
establishing non-constancy. The remaining claims about $m_q$ are immediate. 
\end{proof}
Now we are in a position to use the Bonami-Poornima theorem, which is recalled here.

\begin{lemma}[Theorem 2, \cite{BP87}]
    \label{Bonami Poornima Theorem}
    Let $d \geq 2$ and $k \geq 1$ be integers. Suppose that $m \in C^{\infty}(\R^d \setminus \{0\})$ is continuous and homogeneous of degree zero. Then there exists a constant $C>0$ such that 
$$
\norm{m(D)f}{\dot{W}^{k,1}} \leq C\norm{f}{\dot{W}^{k,1}}, \quad \text{for every } f\in C_c^{\infty}(\R^d),
$$
if and only if $m$ is constant.
\end{lemma}

\begin{remark}
    We remark that the forward direction of the above lemma can be strengthened to only considering $f\in C_c^\infty(\mathbb{R}^d)$ which are also \textit{mean-free}, i.e. $\int_{\mathbb{R}^d}f(y)\,dy =0$. To see this, let $m$ be non-constant and satisfy the assumptions of Lemma \ref{Bonami Poornima Theorem}, and, given $L>0$, let $g\in C^{\infty}_c(\R^d)$ be such that $\| m(D)g \|_{\dot{W}^{k,1}} > L\|g\|_{\dot{W}^{k,1}}$. For fixed $a\in\R^d$ and for each $x\in\R^d$, set 
$$f_a(x) = g(x) - g(x-a),$$ which is mean-free. 
Then for $|\alpha|=k$
$$ D^{\alpha}m(D)f_a(x) = (D^{\alpha}m(D)g)(x) - (D^{\alpha}m(D)g)(x-a).$$ 
Recall that for $h\in L^1(\R^d)$, 
$$\lim_{|a|\rightarrow\infty} \| h - h (\cdot - a)\|_{L^{1}} = 2\| h \|_{L^1},$$ 
so that, 
$$\lim_{|a|\rightarrow\infty} \frac{\|m(D)f_a\|_{\dot{W}^{k,1}}}{\|f_a\|_{\dot{W}^{k,1}}} =  \frac{\| m(D)g \|_{\dot{W}^{k,1}}}{\|g\|_{\dot{W}^{k,1}}} > L.$$
We may therefore choose $|a|$ sufficiently large to ensure that $$\|m(D)f_a\|_{\dot{W}^{k,1}} >L\|f_a\|_{\dot{W}^{k,1}}.$$
This justifies our mean-free choice in the statement of Theorem \ref{Packet Family Theorem} below.
\end{remark}

Now assume $\overline{\theta}_0$ is non-constant, and choose $q_0$ such that $\nabla\overline{\theta}_0(q_0) \neq 0$. Lemmas \ref{multiplier lemma} and \ref{Bonami Poornima Theorem} imply that one can choose $\eta_0\in C^{\infty}_c(\R^2)$ with $\|\eta_0\|_{\dot{W}^{k,1}}\leq 1$ so that the building block $\psi_{q_0}$ generated by (\ref{Packet Equation}) has large $\dot{W}^{k,1}(\R^2)$ norm. In our construction of the perturbation $\theta_0$ of $\overline{\theta}_0$ in Theorem \ref{norm inflation key}, we will require a collection of many building blocks $\psi_q$ with comparably large norm. The following theorem gives this collection by demonstrating continuity of the map $(q,t) \rightarrow \psi_q(t)$ from $\overline{B_\rho(q_0)}\times[0,\tau]$ to $H^{k+1}(\mathbb{R}^2)$ for $\rho$ sufficiently small.         
\begin{theorem}\label{Packet Family Theorem}
    For every $Q>0$ and $\tau>0$ there are mean-free $\eta_0\in C^\infty_c(\mathbb{R}^2)$, $R>0$, $q_0\in B_R$, and $0<\delta<R-|q_0|$, such that, if we define $\psi_q$ as in \eqref{Packet Equation} for $q\in \overline{B_\delta(q_0)}$, then 
    \begin{equation}\label{PFT1}\|\eta_0\|_{\dot{W}^{k,1}}\leq 1,\end{equation}
    and\begin{equation}\label{PFT2}\sum_{|\alpha|=k}\int_{B_R}| D^\alpha\psi_{q}(\tau,y)|\,dy> 4Q \quad \text{ for all } q\in \overline{B_\delta(q_0)}.
    \end{equation}
    Furthermore, for any $\rho>0$, the family $\mathcal{F} = \{\psi_q(t) :q\in \overline{B_\rho(q_0)},0\leq t \leq \tau \}$ is compact in $H^{k+1}(\mathbb{R}^2)$, there exists a constant $C>0$ such that
    \begin{equation}\label{PFT3}
        \sup_{q,t}\Big\{\|\psi_q(t)\|_{H^{k+1}}+\sum_{1\leq |\alpha|\leq k+1}\|\langle y\rangle ^2 D^\alpha \psi_q(t)\|_{L^2}+\sum_{ |\alpha|\leq k}\|\langle y\rangle  D^\alpha R^\perp \psi_q(t)\|_{L^2}\Big\}<C,
    \end{equation}
    and, in particular, for every $L>0$,
    \begin{equation}\label{PFT4}
        \sup_{q,t}\sum_{|\alpha|=k}\int_{|y|\geq L}| D^{\alpha} \psi_q(t,y)|\,dy\leq C/L.
    \end{equation}
    Here the suprema are taken over $(q,t)\in\overline{B_\rho(q_0)}\times [0,\tau]$.
\end{theorem}
\begin{proof}
    Let $Q, \tau>0$ and choose $q_0\in \mathbb{R}^2$ with $\nabla \overline{\theta}_0(q_0)\not=0$. Lemmas \ref{multiplier lemma} and \ref{Bonami Poornima Theorem} give a mean-zero $\eta_0\in C_c^\infty(\mathbb{R}^2)$ such that, for $R>0$ sufficiently large, 
    \begin{equation}\label{choice of eta}
    \|\eta_0\|_{\dot{W}^{k,1}}\leq 1 \quad \text{and} \quad \sum_{|\alpha|=k}\int_{B_R}| D^\alpha\psi_{q_0}(\tau,y)|\,dy\geq 8Q,
    \end{equation}
    which establishes \eqref{PFT1}. To establish \eqref{PFT2}, we fix $\rho>0$ and show that the map $(q,t)\to \psi_q(t)$ from $\overline{B_\rho(q_0)}\times[0,\tau]$ to $H^{k+1}(\mathbb{R}^2)$ is continuous.  Sufficiently small choice of $\delta$ will then give \eqref{PFT2} for our given $\eta_0$.\\ 
    
\noindent {\bf Continuity of the map $(q,t)\rightarrow B_q(t)$.} In view of (\ref{ML:psiq solution formula}), we are motivated to first show continuity of the map $(q,t)\rightarrow B_q(t)$.  Let $(q_n,t_n)\in \overline{B_\rho(q_0)}\times[0,\tau]$ converge to some $(q,t)$. For fixed $s\in [0,\tau]$, by the Mean Value Theorem and smoothness of the background flow $\overline{u}$, 
    $$ |\overline{u}(s,X_{q_n}(s))-\overline{u}(s,X_q(s))| \leq \|\nabla \overline{u}(s) \|_{L^{\infty}}| X_{q_n}(s)-X_q(s) |.$$
    Thus, by (\ref{qtrajectory}) and Gr\"onwall's lemma,
    \begin{equation}
        \label{PFT:trajectory conv q}
        \sup_{0\leq s\leq \tau}|X_{q_n}(s)-X_q(s)|\leq |q_n-q|\exp\left\{\int_0^\tau\|\nabla \overline{u}(s)\|_\infty\,ds\right\}\to 0.
    \end{equation}
    It follows that 
    \begin{equation}
        \label{PFT:trajectory conv qt}
        \begin{split}
            &|X_{q_n}(t_n)-X_q(t)|\leq |X_{q_n}(t_n)-X_q(t_n)|+|X_{q}(t_n)-X_q(t)|
            \\ 
            &\qquad \leq |q_n-q|\exp\left\{\int_0^\tau\|\nabla\overline{u}(s)\|_\infty\,ds\right\}+|t_n-t|\|\overline{u}\|_{L^\infty([0,\tau]\times\mathbb{R}^2)}
            \rightarrow 0.
        \end{split}
    \end{equation}
    Again by smoothness of $\overline{u}$, we have 
    \begin{equation}\label{PFT:nabla background conv qt}
    |\nabla \overline{u}(t_n,X_{q_n}(t_n))-\nabla \overline{u}(t,X_q(t))| + \sup_{0 \leq s \leq \tau} |\nabla\overline{u}(s,X_{q_n}(s))-\nabla\overline{u}(s,X_q(s))| \to 0.
    \end{equation}
By (\ref{definition of Bq}), we have
\begin{equation}
\begin{split}
    &B_{q_n}(s)-B_{q}(s) = \int_0^s\nabla \overline{u}(\sigma,X_{q_n}(\sigma))\left[B_{q_n}(\sigma)-B_{q}(\sigma)\right]\,d\sigma\\
    &\qquad\qquad  +\int_0^s\left[\nabla\overline{u}(\sigma,X_{q_n}(\sigma))-\nabla \overline{u}(\sigma,X_q(\sigma))\right]B_q(\sigma)\,d\sigma.
\end{split}    
\end{equation}
An application of Gr\"onwall's lemma and \eqref{PFT:nabla background conv qt} gives
\begin{equation}
    \label{PFT:Bq conv q}
    \begin{split}
    &\sup_{0\leq s\leq \tau}|B_{q_n}(s)-B_q(s)|\leq \tau e^{\tau\|\nabla \overline{u}\|_{L^\infty_{t,x}}}\sup_{0\leq s\leq \tau}|B_q(s)| \sup_{0\leq s\leq \tau}|\nabla \overline{u}(s,X_{q_n}(s))-\nabla \overline{u}(s,X_q(s))|\\
    &\qquad \leq \tau e^{2\tau\|\nabla \overline{u}\|_{L^\infty_{t,x}}} \sup_{0\leq s\leq \tau}|\nabla \overline{u}(s,X_{q_n}(s))-\nabla \overline{u}(s,X_q(s))|\rightarrow 0,
    \end{split}
\end{equation}
where we applied the bound
$$ \sup_{0\leq s\leq \tau}|B_q(s)| \leq e^{\tau \|\nabla \overline{u}\|_{L^\infty_{t,x}}}. $$
We also have
\begin{equation}
    \label{PFT:Bq conv t} |B_q(t_n)-B_q(t)| = \left|\int_t^{t_n}\nabla \overline{u}(s,X_q(s))B_q(s)\,ds\right|\leq |t_n-t|\sup_{0\leq s\leq \tau}|\nabla \overline{u}(s,X_q(s))B_q(s)|\to0.
\end{equation}
Combining \eqref{PFT:Bq conv q} and \eqref{PFT:Bq conv t} provides $$B_{q_n}(t_n)\to B_{q}(t) \text{ when } (q_n,t_n) \rightarrow (q,t).$$ 
\noindent{\bf Continuity of $(q,t)\to \psi_q(t)$ and proof of (\ref{PFT2}).}
Since $(q,t) \mapsto B_q(t)$ is continuous, from compactness of the set of $\overline{B_\rho(q_0)}\times [0,\tau]$, it follows that 
\begin{equation}
    \label{PFT: Bq bound}
    \sup_{(q,t)\in \overline{B_\rho(q_0)}\times[0,\tau]}\max\left\{|B_q(t)|,|B_q(t)^{-T}|\right\}\leq K
\end{equation}
for some $K\geq 1.$ 

 Now, recall that $$m_q(s,\xi) =\exp\left\{-i\left[\nabla \overline{\theta}_0(q)\cdot \xi^\perp\right]\int_0^s\frac{d\sigma}{|B_q(\sigma)^{-T}\xi|}\right\}.$$
Since the map $(q,t)\rightarrow B_q(t)$ is continuous, it follows that $(q,t)\rightarrow B_q(t)^{-T}$ is continuous, which, in combination with (\ref{PFT: Bq bound}), implies that $(q,t)\rightarrow \frac{1}{|B_q(t)^{-T}|}$ is continuous and bounded below. This continuity, along with smoothness of $\nabla \overline{\theta}_0$, allow us to conclude that $m_q(s,\xi)$ is continuous in $q,s$ for $\xi\not =0$. Since, for each $(q,t)\in \overline{B_\rho(q_0)}\times [0,\tau]$, 
\begin{equation}\label{PFT:Fourierpsiq conv qt}
\begin{split}
&\mathcal{F}_y\big[\psi_{q}(t,y)\big](\xi) = \det \left(B_{q}(t)^{T}\right) m_{q}(t,B_{q}(t)^{T}\xi)\widehat{\eta_0}(B_{q}(t)^{T}\xi)\\
&\qquad\qquad
        =m_{q}(t,B_{q}(t)^{T}\xi)\widehat{\eta_0}(B_{q}(t)^{T}\xi),
\end{split}
\end{equation}
it follows that 
$$ \mathcal{F}_y\big[\psi_{q_n}(t_n,y)\big](\xi) \longrightarrow \mathcal{F}_y\big[\psi_{q}(t,y)\big](\xi) \qquad \text{as }(q_n, t_n) \rightarrow (q,t).$$
Note that \begin{equation}\label{bounds on B}
\|B_{q}(t)^{-T}\| \leq K, \qquad | B_{q}(t)^{T}\xi | \geq \frac{1}{K}|\xi|,
\end{equation}
with $K$ independent of $q,t$ and $\xi$.  Thus by (\ref{PFT:Fourierpsiq conv qt}), the fact that $|m_q(t, \xi)|\equiv 1$, the bounds in \eqref{bounds on B}, and membership of $\widehat{\eta_0}$ to $\mathcal{S}(\R^2)$, for any integer $N>0$, \begin{equation}\label{psihatdecay}
|\widehat{\psi}_q(t,\xi)| =|\widehat{\eta_0}(B_q(t)^T\xi)| \leq C_N(1+|B_q(t)^T\xi|)^{-N}\leq  C_NK^N(1+|\xi|)^{-N}.
\end{equation}
We conclude that
\begin{equation}
    (1+|\xi|^2)^{k+1}|\widehat{\psi_{q_n}}(t_n,\xi)-\widehat{\psi_q}(t,\xi)|^2\in L^1_\xi(\mathbb{R}^2).
\end{equation}
Moreover, as $C_NK^N$ is independent of $q,t$, we can apply the Dominated Convergence Theorem with dominating function
$$ 4C_N^2K^{2N} (1+|\xi|^2)^{k+1}(1+|\xi|)^{-2N},$$
with $N$ sufficiently large, together with Plancherel's theorem to obtain
\begin{equation}
    \label{PFT:psiq Hk+1 convergence}\|\psi_{q_n}(t_n)-\psi_q(t)\|_{H^{k+1}} \to 0,
\end{equation}
which establishes continuity of $(q,t)\to \psi_q(t)$ in $H^{k+1}(\mathbb{R}^2)$ over the set $\overline{B_\rho(q_0)}\times[0,\tau].$ This proves the claim that $\mathcal{F}$ is compact in $H^{k+1}(\mathbb{R}^2)$.

To prove $\eqref{PFT2}$ for $q$ sufficiently close to $q_0$, we apply the Cauchy-Schwarz inequality and \eqref{PFT:psiq Hk+1 convergence} to write
\begin{equation*}
    \label{PFT:proof of PFT2}
    \begin{split}
    &\sum_{|\alpha|=k}\int_{B_R}| D^\alpha\psi_q(\tau,y)|\,dy \geq \sum_{|\alpha|=k}\int_{B_R}| D^\alpha\psi_{q_0}(\tau,y)|\,dy-\sum_{|\alpha|=k}\int_{B_R}| D^\alpha\psi_{q_0}(\tau,y)- D^\alpha\psi_q(\tau,y)|\,dy
    \\  &\qquad\qquad \geq 8Q - C_kR\|\psi_{q_0}(\tau)-\psi_q(\tau)\|_{H^k}
    > 4Q
    \end{split}
\end{equation*}
for $q$ close enough to $q_0$. Define this closeness threshold as $\delta>0$. This proves $\eqref{PFT2}$ for $\eta_0$ as in (\ref{choice of eta}). 
\\ 

\noindent {\bf Proofs of (\ref{PFT3}) and (\ref{PFT4}).} 
From Plancherel's theorem and (\ref{psihatdecay}) it immediately follows that there exists $C>0$ such that $\sup_{q,t}\|\psi_q(t)\|_{H^{k+1}} \leq C$ . Furthermore, since $m_q \in C^{\infty}(\R^d\setminus \{0\})$ is homogeneous of degree zero, (\ref{bounds on B}) and the chain rule give
\begin{equation}\label{derivativembounds}
\begin{split}
&\left| D_{\xi}^{\beta}[m_q(t,B_q(t)^T\xi)]\right| \leq C_{\beta}\|B_q(t)^T\|^{|\beta|}|B_q(t)^T\xi|^{-|\beta|} \leq C_{\beta}|\xi|^{-|\beta|}.
\end{split}
\end{equation}
Since $\eta_0 \in C_c^{\infty}(\R^2)$ is mean-zero, we have that $\widehat{\eta_0}(0)=0$ and, by a Taylor expansion about $0$,
$$
\widehat{\eta_0}(B_q(t)^T\xi) = \widehat{\eta_0}(0) + \nabla\widehat{\eta_0}(0)\cdot B_q(t)^T\xi + o(|B_q(t)^T\xi|) = \mathcal{O}(|\xi|).
$$
Furthermore, every derivative of $\widehat{\eta}_0$ is uniformly bounded. Hence, for $|\alpha|\geq 1$,
\begin{equation}\label{etaderivativebounds}
\left| D_{\xi}^{\alpha} \left(\widehat{\eta}_0(B_q(t)^T\xi)\right)\right| \leq C_{\alpha}|\xi|^{1-|\alpha|}.
\end{equation}
Then, for $0<|\xi| \leq 1$, by (\ref{derivativembounds}), (\ref{etaderivativebounds}), and the Leibniz rule
\begin{equation}\label{psiderivativebounds}
\begin{split}
| D_{\xi}^{\beta}\widehat{\psi_q}(t,\xi)| &= | D_{\xi}^{\beta}[m_q(t,B_q(t)^T\xi)\widehat{\eta_0}(B_q(t)^T\xi)]|\\
&=\left|\sum_{\gamma \leq \beta} C_{\gamma,\beta}  D_{\xi}^{\gamma}\left(m_q(t,B_q(t)^T\xi)\right)  D_{\xi}^{\beta-\gamma}\left(\widehat{\eta_0}(B_q(t)^T\xi)\right)  \right|\\
&\leq \sum_{\gamma \leq \beta} C_{\gamma,\beta} |\xi|^{-|\gamma|} |\xi|^{1-|\beta-\gamma|}\\
&\leq C_{\beta}|\xi|^{1-|\beta|}.
\end{split}
\end{equation}
It follows by another application of the Leibniz rule and (\ref{psiderivativebounds}) that for $1 \leq |\alpha| \leq k+1$ and $0<|\xi|\leq1$,
$$
| D_{\xi}^{\beta}[(i\xi)^{\alpha} \widehat{\psi_q}(t,\xi)]| \leq C_{\alpha,\beta}|\xi|^{|\alpha|+1-|\beta|}.
$$
For $0\leq|\beta|\leq 2$, $|\alpha|+1-|\beta|$ is minimized when $|\alpha|=1$ and $|\beta|=2$, resulting in an exponent of zero. Hence, $ D_{\xi}^{\beta}[(i\xi)^{\alpha} \widehat{\psi_q}(t,\xi)]$ is square-integrable near the origin, and decays rapidly at infinity. It follows from Plancherel's theorem that 
$$
\norm{y^{\beta}D^{\alpha}\psi_q(t)}{L^2} \leq C_{\eta_0}, \quad 1\leq |\alpha| \leq k+1, 0\leq |\beta|\leq 2.
$$
From the Leibniz rule, $\langle y\rangle^2\nabla\psi_q(t) \in H^k(\R^2)$. An identical argument gives $\langle y\rangle R^{\perp}\psi_q(t) \in H^k(\R^2)$ by noting that 
$$
\left| D_{\xi}^{\beta}\left[(i\xi)^{\alpha}\frac{\xi^{\perp}}{|\xi|}\widehat{\psi_q}(t,\xi)\right]\right| \leq C_{\alpha,\beta}|\xi|^{|\alpha|+1-|\beta|}, \quad \text{for } 0\leq |\alpha| \leq k, 0 \leq |\beta|\leq1
$$
establishing \eqref{PFT3}.

Finally, to prove \eqref{PFT4}, we apply the Cauchy-Schwarz inequality and \eqref{PFT3}, which gives
\begin{equation}
    \label{PFT:proof of PFT4}
    \begin{split}
        \sum_{|\alpha|=k}\int_{|y|\geq L}| D^\alpha\psi_q(t,y)|\,dy \leq \left(\int_{|y|\geq L}\langle y\rangle ^{-4}\,dy\right)^{1/2}\sum_{|\alpha|=k}\|\langle y\rangle ^2 D^\alpha\psi_q\|_{L^2}\leq C/L.
    \end{split}
\end{equation}
\end{proof}

\begin{remark}
    The purpose of the weighted norms in \eqref{PFT3} will become clear in the next subsection. In the approximation of the background drift and scalar, the Taylor estimate induces errors of the form $\langle y\rangle^2$ and $\langle y\rangle$, which we 
    control using \eqref{PFT3}.
\end{remark}

\subsection{Tracking Many Blocks}
Fix $Q,\tau,\eta_0,R,\delta,q_0$ from Theorem \ref{Packet Family Theorem}. Choose $J\geq 2$ distinct points $q_1,...,q_J$ in $B_\delta(q_0).$ For $r>0$ to be chosen small and specified later, set
\begin{equation}
    \Psi(t,y) = \sum_{j=1}^J\psi_{q_j}\left(t,y-\frac{X_{q_j}(t)}{r}\right).
\end{equation}
In what follows, $\Psi$ will serve as an approximation to the actual perturbation.  To clarify, fix $a>0$ and consider the quantity
\begin{equation}\label{vdef}
v(t,y) = \frac{\theta(t,ry)-\overline{\theta}(t,ry)}{ar^{k-2}}.
\end{equation}
One can quickly check that the evolution of $v$ is given by
\begin{equation}\label{vevolution1}
\begin{cases}
\partial_tv+R^\perp v\cdot \nabla\overline{\theta}(t,ry)+\frac{\overline{u}(t,ry)}{r}\cdot\nabla v +ar^{k-3}R^\perp v\cdot\nabla v = 0\\
v|_{t=0} = \Psi(0,\cdot). 
\end{cases}
\end{equation}
Moreover, by the chain rule and a change of variables,
$$\|\theta(t)-\overline{\theta}(t)\|_{W^{k,1}}\geq \|(\theta(t)-\overline{\theta}(t))\|_{\dot{W}^{k,1}} = a\| v (t)\|_{\dot{W}^{k,1}}.$$
While we must eventually choose $a$ small, we are motivated by the above to construct $\theta_0$ in such a way that $\| v(t) \|_{\dot{W}^{k,1}}$ becomes large for small $t$. To achieve this goal, we will approximate $v$ by $\Psi$, which we can make large using Theorem \ref{Packet Family Theorem}. 

In Lemma \ref{Lemma Moving Family of Packets} below, we obtain the necessary estimates to control the error between $\Psi$ and the actual evolved perturbation $v$. We also obtain the desired $\dot{W}^{k,1}$ norm inflation for $\Psi$, which will in turn yield a similar inequality for $v$ when $J$ is chosen sufficiently large.

\begin{lemma}
    \label{Lemma Moving Family of Packets}
    There is a constant $C>0$ independent of $J,r,$ and the choice of points $q_j$, such that for all sufficiently small $r>0$ ,
    \begin{equation}\label{MFP1}
        \begin{split}
            &\sup_{t\leq \tau}\|\Psi(t)\|_{H^{k+1}}\leq C\sqrt{J},
\\            &\sup_{t\leq \tau}\Big\|\partial_t\Psi+\frac{\overline{u}(t,r\cdot)}{r}\cdot\nabla\Psi+\nabla\overline{\theta}(t,r\cdot)\cdot R^\perp\Psi\Big\|_{H^k}\leq CJr.
        \end{split}
    \end{equation}
    Furthermore, at time $\tau$, the balls $B_j:=B_R\left(\frac{X_{q_j}(\tau)}{r}\right)$ are disjoint and $\Psi$ satisfies    \begin{equation}\label{MFP2}
       \sum_{|\alpha|=k} \int_{\cup_{j=1}^J B_j}| D^\alpha\Psi(\tau,y)|\,dy\geq 2JQ.
    \end{equation}
    The smallness threshold for $r$ depends on $J$ and the choice of $q_j$.
\end{lemma}
\begin{proof}
{\bf Proof of \eqref{MFP1}$_1$.}
    Since $q_1,...,q_J$ are distinct, their background trajectories remain uniformly separated. That is, for $1\leq j<i\leq J,$
    $$|X_{q_j}(t)-X_{q_i}(t)|\geq |q_j-q_i|\exp\left\{-\tau\|\nabla\overline{u}\|_{L^\infty([0,\tau]\times \mathbb{R}^2)}\right\}>d_J,$$ for some $d_J>0$ depending only on the choice of points $q_j.$ Thus, for $r>0$ sufficiently small, $B_j\cap B_i=\emptyset$ for $j\not=i$.
    
    By Theorem \ref{Packet Family Theorem}, the family $\{\psi_{q_j}(t):1\leq j\leq J,0\leq t\leq \tau\}$ is compact in $H^{k+1}(\mathbb{R}^2)$. Thus, the equitightness condition of the Fr\'echet-Kolmogorov-Riesz theorem provides uniformly small $L^2(\mathbb{R}^2)$ tails of all derivatives of $\psi_{q_j}$ up to order $k+1$; that is,
    \begin{equation}\label{MFP: equitightness}
        \lim_{L\to \infty}\sup_{0\leq t\leq \tau}\sup_{1\leq j\leq J}\sup_{|\alpha|\leq k+1}\int_{|y|\geq L}| D^\alpha\psi_{q_j}(t,y)|^2\,dy = 0.
    \end{equation}
    Now, taking the $H^{k+1}(\mathbb{R}^2)$ inner product of two terms of the sum defining $\Psi$, for fixed $L>0$ we find that 
    \begin{equation*}
    \begin{split}
        &\qquad \Bigg\langle\psi_{q_j}\left(t,y-\frac{X_{q_j}(t)}{r}\right),\psi_{q_i}\left(t,y-\frac{X_{q_i}(t)}{r}\right)\Bigg\rangle_{H^{k+1}} \\ &=\sum_{|\alpha|\leq k+1} \bigg [\int_{B_{L}\left(\frac{X_{q_j}(t)}{r}\right)} D^\alpha\psi_{q_j}\left(t,y-\frac{X_{q_j}(t)}{r}\right) D^\alpha\psi_{q_i}\left(t,y-\frac{X_{q_i}(t)}{r}\right)\,dy \\
        &+\int_{\mathbb{R}^2\backslash B_{L}\left(\frac{X_{q_j}(t)}{r}\right)} D^\alpha\psi_{q_j}\bigg(t,y-\frac{X_{q_j}(t)}{r}\bigg) D^\alpha\psi_{q_i}\left(t,y-\frac{X_{q_i}(t)}{r}\right)\,dy \bigg]
        \\ &\leq\sum_{|\alpha|\leq k+1}\bigg (\| D^\alpha \psi_{q_j}(t)\|_{L^2}\bigg(\int_{|z|\leq L}\left| D^\alpha\psi_{q_i}\left(t,z+\frac{X_{q_j}(t)-X_{q_i}(t)}{r}\right)\right|^2\,dz\bigg)^{1/2}
        \\ &\qquad \qquad\qquad +\left(\int_{|z|\geq L}| D^\alpha\psi_{q_j}(t,z)|^2\,dz\right)^{\frac{1}{2}}\left\| D^\alpha\psi_{q_i}(t)\right\|_{L^2} \bigg ),
        \end{split}
    \end{equation*}
    where we have applied the Cauchy-Schwarz inequality and the change of variables $z=y-\frac{X_{q_j}(t)}{r}$. We bound the $L^2$ norms appearing above using \eqref{PFT3}. To make the remaining integrals small, we use \eqref{MFP: equitightness} to first choose $L$ large, and then, given this $L$, we choose $r$ small. More explicitly, let $\varepsilon>0$, and choose $r>0$ sufficiently small so that
    \begin{equation}
    \begin{split}
        \label{MFP:Proof of MFP1,1}
        \|\Psi(t)\|_{H^{k+1}}^2 &= \sum_{j=1}^J \|\psi_{q_j}(t)\|_{H^{k+1}}^2+\sum_{i\not= j}\Bigg\langle\psi_{q_j}\left(t,y-\frac{X_{q_j}(t)}{r}\right),\psi_{q_i}\left(t,y-\frac{X_{q_i}(t)}{r}\right)\Bigg\rangle_{H^{k+1}} 
        \\ &\leq C^2J+\varepsilon
        \end{split}
    \end{equation}
    where $C>0$ is the constant in the statement of Theorem \ref{Packet Family Theorem}. Thus, there is a constant, which we relabel $C>0$, for which 
    $$\sup_{0\leq t\leq \tau}\|\Psi(t)\|_{H^{k+1}}\leq C\sqrt{J}, $$ which is \eqref{MFP1}$_1$.\\
\\
    \noindent {\bf Proof of \eqref{MFP1}$_2$.} Fix $j$ satisfying $1\leq j \leq J$. To simplify notation, in what follows we set $q=q_j$. As above, we introduce the translated variable $z=y-X_{q}(t)/r$. Since $\psi_q$ satisfies \eqref{Packet Equation}, given by
    $$
    \partial_t\psi_q + (\nabla \overline{u}(t,X_q(t))y)\cdot \nabla \psi_q + R^{\perp}\psi_q\cdot \nabla \overline{\theta}(t,X_q(t))=0,
    $$
    the chain rule gives 
    \begin{equation*}
    \begin{split}
    \label{no substitute}
        & \partial_t\left(\psi_q(t,y-X_{q}(t)/r)\right)+\frac{\overline{u}(t,ry)}{r}\cdot\nabla\psi_q(t,y-X_q(t)/r)+\nabla\overline{\theta}(t,ry)\cdot R^\perp\psi_q(t,y-X_{q}(t)/r) \\&= \partial_t\psi_q(t,z)-\frac{X'_q(t)}{r}\cdot\nabla\psi_q(t,z)+\frac{\overline{u}(t,X_q(t)+rz)}{r}\cdot\nabla \psi_q(t,z)+\nabla\overline{\theta}(t,X_q(t)+rz)\cdot R^\perp\psi_q(t,z)
        \\ & = -\nabla\overline{u}(t,X_q(t))z\cdot \nabla\psi_q(t,z)+\frac{\overline{u}(t,X_q(t)+rz)-\overline{u}(t,X_q(t))}{r}\cdot\nabla \psi_q(t,z)
        \\ &\qquad \qquad+\big(\nabla\overline{\theta}(t,X_q(t)+rz)-\nabla\overline{\theta}(t,X_q(t))\big)\cdot R^\perp\psi_q(t,z)
        \\ 
        & =\Xi_{q,r}(t,z)\cdot\nabla\psi_q(t,z) + \Theta_{q,r}(t,z)\cdot R^\perp \psi_q(t,z),
        \end{split}
    \end{equation*}
    where 
    \begin{equation*}
    \begin{split}
    &\Xi_{q,r}(t,z) = \frac{\overline{u}(t,X_q(t)+rz)-\overline{u}(t,X_q(t))}{r}-\nabla\overline{u}(t,X_q(t))z, \\
    &\Theta_{q,r}(t,z) = \nabla\overline{\theta}(t,X_q(t)+rz)-\nabla\overline{\theta}(t,X_q(t)).
    \end{split}
    \end{equation*}
    By Taylor's theorem, there exists $C_k>0$, depending on $k$ and $\overline{\theta}_0$ but not on $q,r,$ or $z,$ such that for every $|\beta|\leq k$,
    \begin{equation*}
    \label{taylors theorem}
        \begin{split}
            &\left | D_z^\beta\Xi_{q,r}(t,z)\right|\leq C_k r\langle z\rangle^2, \qquad \left| D_z^\beta\Theta_{q,r}(t,z)\right|\leq C_kr\langle z\rangle.
        \end{split}
    \end{equation*}
     Now if we sum over $q_j$ for $ j=1,...,J$ and apply the $H^k(\mathbb{R}^2)$ norm, we obtain
    \begin{equation*}\begin{split}
        &\|\partial_t\Psi+\frac{\overline{u}(t,r\cdot)}{r}\cdot\nabla \Psi+\nabla \overline{\theta}(t,r\cdot)\cdot R^\perp \Psi\|_{H^k}
        \\ &\leq \sum_{j=1}^J\bigg\|\partial_t\left[\psi_{q_j}\left(t,\cdot-\frac{X_{q_j}(t)}{r}\right)\right]+\frac{\overline{u}(t,r\cdot)}{r}\cdot\nabla \psi_{q_j}\left(t,\cdot-\frac{X_{q_j}(t)}{r}\right)\\
        &\qquad\qquad  +\nabla \overline{\theta}(t,r\cdot)\cdot R^\perp \psi_{q_j}\left(t,\cdot-\frac{X_{q_j}(t)}{r}\right)\bigg\|_{H^k}
        \\ &\leq \sum_{j=1}^JC_k\left\|r\left\langle \cdot-\frac{X_{q_j}(t)}{r}\right\rangle^2\nabla \psi_{q_j}\left(t,\cdot-\frac{X_{q_j}(t)}{r}\right)\right\|_{H^k} \\
        &\qquad \qquad  + \sum_{j=1}^JC_k\left\|r\left\langle \cdot-\frac{X_{q_j}(t)}{r}\right\rangle R^\perp \psi_{q_j}\left(t,\cdot-\frac{X_{q_j}(t)}{r}\right)\right\|_{H^k} 
        \\ &\leq CC_k Jr,
        \end{split}
    \end{equation*}
    where we used \eqref{PFT3} and translation invariance of the $H^k(\mathbb{R}^2)$ norm. This proves assertion \eqref{MFP1}$_2$.
    \\
\\
    \noindent {\bf Proof of \eqref{MFP2}.} By \eqref{PFT4}, we have
    \begin{equation}\label{MFP:far term}
       \sum_{|\alpha|=k}\sum_{j=1}^J \sum_{j\not=i}\int_{|x|\geq d_J/r-R}| D^\alpha\psi_{q_j}(\tau,x)|\,dx \leq C_kJ(J-1)\cdot\frac{C}{d_J/r-R}.
    \end{equation}
    For fixed $J$, choose $r$ small enough to ensure that the expression on the right side of \eqref{MFP:far term} is at most $2JQ$. Then for $|\alpha|=k$, by the reverse triangle inequality, for each $j=1,...,J$,
    $$| D^\alpha\Psi(\tau,y)|\geq\left| D^\alpha\psi_{q_j}\left(\tau,y-\frac{X_{q_j}(\tau)}{r}\right)\right|-\sum_{i\not=j}\left| D^\alpha\psi_{q_i}\left(\tau,y-\frac{X_{q_i}(\tau)}{r}\right)\right|.$$
    Integrating over $\bigcup_{j=1}^JB_j$ and using the disjointness of the $B_j$, we obtain
    \begin{equation*}
        \label{MFP: setup proof of MFP2}
        \begin{split}
            &\int_{\cup_{j=1}^JB_j} | D^\alpha\Psi(\tau,y)|\,dy = \sum_{j=1}^J\int_{B_j}| D^\alpha\Psi(\tau,y)|\,dy
            \\ &\qquad\geq \sum_{j=1}^J\int_{B_j}\left| D^\alpha \psi_{q_j}\left(\tau,y-\frac{X_{q_j}(\tau)}{r}\right)\right|\,dy
            -\sum_{j=1}^J\sum_{i\not=j}\int_{B_j}\left| D^\alpha\psi_{q_i}\left(\tau,y-\frac{X_{q_i}(\tau)}{r}\right)\right|\,dy\\
              & \qquad= \sum_{j=1}^J\int_{B_R(0)}| D^\alpha\psi_{q_j}(\tau,z)|\,dz
             -\sum_{j=1}^J\sum_{i\not=j}\int_{B_R\left(\frac{X_{q_j}(\tau)-X_{q_i}(\tau)}{r}\right)}| D^\alpha\psi_{q_i}(\tau,x)|\,dx,
        \end{split}
    \end{equation*}
    where in the final equality we made the change of variables $z=y-\frac{X_{q_j}(\tau)}{r}$ in the first integral, and $x=y-\frac{X_{q_i}(\tau)}{r}$ in the latter. Note that if we choose $r>0$ sufficiently small so that $r<d_J/(2R)$, then $x\in B_{R}\left(\frac{X_{q_j}(\tau)-X_{q_i}(\tau)}{r}\right)$ implies $|x|\geq d_J/r-R$. Thus, summing over $|\alpha|=k$ and applying \eqref{MFP:far term} and \eqref{PFT2} yields
    \begin{equation*}
        \sum_{|\alpha|=k}\int_{\cup_{j=1}^JB_j}| D^\alpha \Psi(\tau,y)|\,dy\geq 4QJ - 2QJ = 2QJ,
    \end{equation*}
    which is \eqref{MFP2}.
\end{proof}

\subsection{Choosing the Initial Data}
Now for $a,r>0$, we define the initial data \begin{equation}
    \label{initial data definition}
    \theta_0(x)=\overline{\theta}_0(x)+ar^{k-2}\sum_{j=1}^J\eta_0\left(\frac{x-q_j}{r}\right). 
\end{equation}

Having established the desired norm inflation estimate on $\Psi$ in Lemma \ref{Lemma Moving Family of Packets}, we now show that the exact perturbation $v$, as defined in (\ref{vdef}), can be adequately approximated by $\Psi$ in the $H^k$ norm. This approximation will allow us to upgrade the lower bound \eqref{MFP2} from $k$ derivatives of $\Psi$ to $k$ derivatives of $v$ in $L^1(\mathbb{R}^2)$.  

\begin{lemma}\label{Error Lemma}
    Fix $b>0$ and set $a=b/J$. Assume $\theta_0$ and $\bar{\theta}_0$ are as in (\ref{initial data definition}), and assume $\theta$  and $\bar{\theta}$ are the short time solutions to (SQG) generated from $\theta_0$ and $\bar{\theta}_0$, respectively. Let $v$ be as in (\ref{vdef}).
    Then for all sufficiently large $J$ and all sufficiently small $r>0$ as given by Lemma \ref{Lemma Moving Family of Packets}, \begin{equation}\label{EL1}
            \sup_{0\leq t\leq \tau}\|v(t)\|_{H^k}\leq C\sqrt{J},
            \end{equation}
            \begin{equation}\label{EL2}
         \sup_{0\leq t\leq \tau}\|v(t)-\Psi(t)\|_{H^k}\leq CJ(a+r),
           \end{equation}
    for $C>0$ depending on $\eta_0,\overline{\theta}_0,b,k,\tau$, but independent of $J,r$.
\end{lemma}
\begin{remark}
Despite the requirement for large $J$ in Lemma \ref{Error Lemma}, when applying (\ref{EL2}) to the proof of Theorem \ref{norm inflation key}, we will choose $a$ and $r$ sufficiently small to ensure that $\Psi$ adequately approximates $v$.  
\end{remark}
\begin{proof}[Proof of Lemma \ref{Error Lemma}]
First recall that the evolution of $v$ is given as in (\ref{vevolution1}); that is, 
\begin{equation}\label{vevolution}
\begin{cases}
\partial_tv+R^\perp v\cdot \nabla\overline{\theta}(t,ry)+\frac{\overline{u}(t,ry)}{r}\cdot\nabla v +ar^{k-3}R^\perp v\cdot\nabla v = 0\\
v|_{t=0} = v_0, 
\end{cases}
\end{equation}
with $v_0 = \Psi(0,\cdot).$ By a standard energy argument in which we apply $D^{\alpha}$, $|\alpha|\leq k$, to \eqref{vevolution} and test against $D^\alpha v$, we obtain (see Lemma \ref{energy estimate 1})
   \begin{equation}\label{EL:to be proved in appendix}
       \|v(t)\|_{H^k}\leq \frac{e^{Ct}\|v_0\|_{H^k}}{1-Ca\|v_0\|_{H^k}\int_0^te^{Cs}\,ds },
   \end{equation}
    where $C>0$ is independent of $a,r,$ and $v_0$. Since $v_0(y) = \sum_{j=1}^J\eta_0\left(y-\frac{q_j}{r}\right)$, by taking $r>0$ sufficiently small so that the supports of the summands are disjoint,
   $$\|v_0\|_{H^k}^2= \sum_{j=1}^J\|\eta_0\|_{H^k}^2 = J\|\eta_0\|_{H^k}^2.$$ Since $a = b/J$ with $b>0$ fixed, $a\|v_0\|_{H^k}=\|\eta_0\|_{H^k}b/\sqrt{J}$. Choose $J$ sufficiently large so that $$Ca\|v_0\|_{H^k}\int_0^\tau e^{Cs}\,ds=\frac{Cb\|v_0\|_{H^k}}{J}\int_0^\tau e^{Cs}\,ds<\frac{1}{2}.$$ Substituting this into \eqref{EL:to be proved in appendix} establishes \eqref{EL1}.

   To prove \eqref{EL2}, for ease of notation we set:
   \begin{equation*}
       \begin{split}
           U_r(t,y)  = \frac{\overline{u}(t,ry)}{r},\quad G_r (t,y)= \nabla\overline{\theta}(t,ry),\quad w = v-\Psi
       \end{split}
   \end{equation*}
   and 
   \begin{equation*}
       E = \partial_t\Psi+U_r\cdot\nabla \Psi +G_r\cdot R^\perp\Psi.
   \end{equation*}
   Subtracting $E$ from \eqref{vevolution}, we obtain
   \begin{equation}
   \label{w evolution}
       \partial_tw+(U_r+ar^{k-3}R^\perp v)\cdot\nabla w +G_r\cdot R^\perp w +ar^{k-3}R^\perp w\cdot \nabla \Psi = -E-ar^{k-3}R^\perp \Psi\cdot\nabla\Psi
   \end{equation}
   with $w|_{t=0} \equiv 0$. An application of Lemma \ref{energy estimate 2} gives
   \begin{equation}
       \label{EL: error energy estimate appendix}
       \begin{split}
       \frac{d}{dt}\|w\|_{H^k}\leq C\big[1&+ar^{k-3}\|v\|_{H^k}+ar^{k-3}\|\Psi\|_{H^{k+1}}\big]\|w\|_{H^k}
       \\ & +\|E\|_{H^k}+Car^{k-3}\|\Psi\|_{H^k}\|\Psi\|_{H^{k+1}},
       \end{split}
   \end{equation}
   where $C>0$ is independent of $a,r,$ and $v_0$. Since $a=b/J$, by \eqref{MFP1}$_1$ and \eqref{EL1}, the coefficient on $\|w\|_{H^k}$ can be bounded as follows,
   \begin{equation*}
       1+ar^{k-3}\|v\|_{_{H^k}}+ar^{k-3}\|\Psi\|_{H^{k+1}} \leq 1+Car^{k-3}\sqrt{J}= 1+Cr^{k-3}b/\sqrt{J},
   \end{equation*}
   which is less than $2$ for sufficiently large $J$, uniformly in $r \in (0,1]$. Moreover, by \eqref{MFP1}$_2$
   \begin{equation}\label{EL:E est}\|E\|_{H^k}\leq CJr,
   \end{equation}
   and by \eqref{MFP1}$_1$, again, we have
   \begin{equation}\label{EL:pert}ar^{k-3}\|\Psi\|_{H^k}\|\Psi\|_{H^{k+1}}\leq CaJ.
   \end{equation}
   Substituting the above estimates into (\ref{EL: error energy estimate appendix}), applying Gr\"onwall's lemma, and using that $w|_{t=0} \equiv 0$, we conclude that $\|v(t)-\Psi(t)\|_{H^k}\leq CJ(a+r)$, proving \eqref{EL2}.
\end{proof}
We are ready to prove the final lemma, which provides the $L^1(\mathbb{R}^2)$ lower bound on the $k^{th}$ order derivatives of $\theta(\tau)-\overline{\theta}(\tau)$ in terms of $b,Q,J,$ and $r$. This is the last ingredient needed for the proof of Theorem \ref{norm inflation key}.
\begin{lemma}
    \label{L1 lower bounds lemma}
    With all the same notation and choices as in the previous sections, 
    \begin{equation}\label{LB1}
        \|\theta(\tau)-\overline{\theta}(\tau)\|_{W^{k,1}}\geq 2bQ-C_R\left(\frac{b^2}{\sqrt{J}}+b\sqrt{J}r\right) 
    \end{equation}
    where $C_R>0$ is a constant which depends only on $\eta_0,\overline{\theta}_0, b,k,\tau$, but is independent of $J$ and $r$. 
\end{lemma}
\begin{proof}
    By the Cauchy-Schwarz inequality and Lemma \ref{Error Lemma}, we have 
    \begin{equation*}
    \begin{split}
    \sum_{|\alpha|=k}\int_{\cup_{j=1}^J B_j}|D^{\alpha}(v-\Psi)(\tau,y)|\,dy\leq \Big|\bigcup_{j=1}^J B_j\Big|^{1/2} \|v(\tau)-\Psi(\tau)\|_{H^k}\leq C_RJ^{3/2}(a+r).
    \end{split}
    \end{equation*}
    An application of \eqref{MFP2} and the reverse triangle inequality now gives
    \begin{equation}\label{vlowerbound}\sum_{|\alpha|=k}\int_{\cup_{j=1}^J B_j}|D^{\alpha} v(\tau)|\,dy\geq 2JQ - C_{R}J^{3/2}(a+r).
    \end{equation}
    By definition, $\theta(\tau,ry)-\overline{\theta}(\tau,ry) = ar^{k-2}v(\tau,y)$, so that by the chain rule and a change of variables,
    $$\|\theta(\tau)-\overline{\theta}(\tau)\|_{W^{k,1}}\geq \|(\theta(\tau)-\overline{\theta}(\tau))\|_{\dot{W}^{k,1}} = a\| v\|_{\dot{W}^{k,1}}.$$ Thus, applying (\ref{vlowerbound}), we conclude that $$\|\theta(\tau)-\overline{\theta}(\tau)\|_{W^{k,1}} \geq 2aJQ-C_RaJ^{3/2}(a+r).$$ Using $aJ = b$ gives \eqref{LB1}.
\end{proof}
\begin{remark}\label{Remark dilution}
    To clarify the roles of the different parameters in the estimate \eqref{LB1}, we point out the contributions made to each term on the right-hand side. The positive term $2bQ$ comes directly from the lower bound on $\Psi$ in \eqref{MFP2}. The defect terms come from the errors induced when approximating the actual evolved perturbation $v$ with $\Psi$. In particular, the term $b^2/\sqrt{J}$ arises due to the estimate \eqref{EL:pert} on the nonlinear term appearing in the evolution of the error $v-\Psi$ \eqref{w evolution}. The rest of the error is contained in the term $b\sqrt{J}r$. The first defect term is precisely where we see the mechanism of \textit{many-packet dilution} referenced in the introduction. As we place more packets in the perturbation, the effect of the nonlinear term has size on the order of $1/\sqrt{J}$, which vanishes as $J$ is made large. On the other hand, increasing $J$ increases the remaining error $b\sqrt{J}r$. We counteract this below by choosing the spatial scale $r>0$ sufficiently small after fixing $J$.
\end{remark}

\subsection{Proof of Theorem \ref{norm inflation key}}
 We now conclude the proof of Theorem \ref{norm inflation key}, which amounts to selecting all parameters in the correct order to prevent circularity and ensuring that the perturbed initial data $\theta_0$ remains supported in the same ball $B_L(0)$ as $\overline{\theta}_0$. Let $\varepsilon,M>0$ and fix $L>0$ such that 
$$
\supp \overline{\theta_0} \subseteq B_L(0).
$$
Choose $$0<b<\varepsilon/2,\qquad Q>M/b.$$
Use Theorem \ref{Packet Family Theorem} to choose $\eta_0,R,q_0, $ and $\delta$. Since $\nabla \overline{\theta}_0(q_0) \neq 0$ we have $q_0 \in B_L(0)$. Set
$$
d=L-|q_0|>0.
$$
Shrinking $\delta$ if necessary, we may assume that 
$$
0<\delta<d/2.
$$
Note that all conclusions of Theorem \ref{Packet Family Theorem} remain valid for this smaller value of $\delta$. Also, choose $R_0>0$ so that 
$$
\supp \eta_0 \subseteq B_{R_0}(0).
$$
Choose $J>0$ large enough to apply Lemma \ref{Error Lemma} and to satisfy
$$C_Rb^2/\sqrt{J}\leq bQ/2,$$
where $C_R>0$ is the constant from Lemma \ref{L1 lower bounds lemma}. Select $J$ distinct points within $B_\delta(q_0)$. Now choose $r>0$ sufficiently small to apply Lemmas \ref{Lemma Moving Family of Packets} and \ref{Error Lemma} and to guarantee
$$C_Rb\sqrt{J}r\leq bQ/2,\qquad \sum_{|\alpha|\leq k-1}r^{k-|\alpha|}\| D^\alpha \eta_0\|_{L^1}\leq 1, \qquad rR_0<\frac{d}{2}.$$
Finally, let $a=b/J$ and define $\theta_0$ by 
$$
\theta_0(x)=\overline{\theta}_0(x)+ar^{k-2}\sum_{j=1}^J\eta_0\left(\frac{x-q_j}{r}\right).
$$For each $j$,
$$
\supp\left[\eta_0\left(\frac{\cdot-q_j}{r}\right)\right]=q_j+r\supp\eta_0 \subseteq B_{\delta+rR_0}(q_0) \subseteq B_L(0),
$$
where in the last set containment we used that $\delta+rR_0<d=L-|q_0|$. It follows that $\supp \theta_0 \subseteq B_L(0)$.
Now, since $\|\eta_0\|_{\dot{W}^{k,1}}\leq 1$, 
$$\|\theta_0-\overline{\theta}_0\|_{W^{k,1}}\leq aJ\sum_{|\alpha|\leq k}r^{k-|\alpha|}\| D^\alpha\eta_0\|_{L^1}\leq 2b <\varepsilon.$$
Lemma \ref{Error Lemma} supplies the existence of a smooth solution $\theta$ through time $\tau$. Lemma \ref{L1 lower bounds lemma} gives
\begin{equation}
    \|\theta(\tau)-\overline{\theta}(\tau)\|_{W^{k,1}}\geq 2bQ-C_R\frac{b^2}{\sqrt{J}} - C_Rb\sqrt{J}r\geq bQ>M,
\end{equation}
which concludes the proof.
\section{Proof of Theorem \ref{BLOW UP}: Instantaneous Blow-Up}
\label{instant blow up section}
We now use Theorem \ref{norm inflation key} to construct a single initial datum whose solution exhibits instantaneous blow-up of its $W^{k,1}(\R^2)$ norm. The main difficulty in this construction is to preserve the large norms obtained at earlier stages while introducing additional perturbations. To achieve this, we first establish a common existence time interval and $L^2(\R^2)$-stability estimate for solutions whose initial data are sufficiently close to the background scalar in $W^{k,1}(\R^2)$. We then detect large $W^{k,1}(\R^2)$ norms using mollifying and localizing functionals, which depend continuously on the solution in $L^2(\R^2)$. This allows us to construct a convergent sequence of smooth initial data with increasingly large norms at times tending to zero, while ensuring that these lower bounds persist for the limiting solution.
\subsection{Generating a Sequence of Initial Data}
In what follows, we fix $k \geq 3$ and let $\eps>0$ and $\overline{\theta}_0 \in C_c^{\infty}(\R^2)$ be given. For $r>0$ satisfying 
$$r<\min\left\{\eps,\frac{1}{2}\norm{\overline{\theta}_0}{W^{k,1}}\right\},$$
we consider the family 
$$
\mathcal{G}= \left\{f \in W^{k,1}(\R^2) : \norm{f-\overline{\theta_0}}{W^{k,1}}<r\right\}.
$$
Note that the assumption on $r$ above ensures that the functions in $\mathcal{G}$ are non-constant. 

\begin{lemma}
\label{H21 LWP and L2 stab}
There exists a single $T>0$ such that every $f \in \mathcal{G}$ generates a weak solution (constructed in Appendix \ref{LWP}) $\theta_f$ to \eqref{SQG} on $[0,T]\times \R^2$. Furthermore, there exist constants $C_1$ and $K$ such that 
\begin{equation}
\label{uniform H21 bound ibu}
\sup_{f \in \mathcal{G}} \sup_{0 \leq t \leq T} \norm{\theta_f(t)}{W^{2,(2,1)}} \leq K,
\end{equation}
and, for all $f,g\in \mathcal{G}$, 
\begin{equation}
\label{L2 stab}
\sup_{0 \leq t \leq T} \norm{\theta_f(t)-\theta_g(t)}{L^2} \leq C_1\norm{f-g}{L^2}.
\end{equation}
The constants $C_1$ and $K$ only depend on $\|\overline{\theta}_0\|_{W^{k,1}}$. The space $W^{2,(2,1)}(\R^2)$ is defined in Definition \ref{Lorentz space definition}.
\end{lemma}
\begin{proof}
Let $f \in \mathcal{G}$. From Theorem \ref{lorentz embedding 1}, $W^{k,1}(\R^2) \hookrightarrow W^{2,(2,1)}(\R^2)$. It follows that 
$$
\norm{f}{W^{2,(2,1)}} \leq \norm{\overline{\theta}_0}{W^{2,(2,1)}} + \norm{f-\overline{\theta}_0}{W^{2,(2,1)}} \leq C \left(\norm{\overline{\theta}_0}{W^{k,1}}+r\right).
$$
From Theorem \ref{LWP H21}, there exists $C,T>0$ such that each $f \in \mathcal{G}$ generates a solution $\theta_f$ satisfying 
$$
\sup_{0 \leq t \leq T} \norm{\theta_f(t)}{W^{2,(2,1)}} \leq 2C(\|\overline{\theta}_0\|_{W^{2,(2,1)}}+r),
$$
where $T \leq \frac{1}{C\left(\norm{\overline{\theta}_0}{W^{k,1}}+r\right)}$, establishing \eqref{uniform H21 bound ibu} with $K=2C(\|\overline{\theta}_0\|_{W^{2,(2,1)}}+r)
$. 

Now, given $f,g \in \mathcal{G}$, we consider the difference of their solutions, $w=\theta_f-\theta_g$ which satisfies
\begin{equation}
\label{H21 l2 stab eq}
\begin{cases}
\partial_tw + u_f\cdot \nabla w + R^{\perp}w \cdot \nabla \theta_g =0,\\
w\vert_{t=0}=f-g.
\end{cases}
\end{equation}
Taking the $L^2$ inner-product of \eqref{H21 l2 stab eq} with $w$, and using that $u_f$ is divergence free yields 
$$
\frac{1}{2}\frac{d}{dt} \norm{w(t)}{L^2}^2 =- \int_{\R^2} (R^{\perp}w(t,x)\cdot \nabla \theta_g(t,x))w(t,x)\,dx \leq C\norm{w(t)}{L^2}^2\norm{\nabla \theta_g(t)}{L^{\infty}}.
$$
From Gr\"{o}nwall's Lemma, setting $w=\theta_f-\theta_g$, it follows that
$$
\norm{\theta_f(t)-\theta_g(t)}{L^2}\leq \norm{f-g}{L^2}\exp\left(C\int_0^t \norm{\nabla \theta_g(\tau)}{L^{\infty}}\,d\tau\right).
$$
By Theorem \ref{lorentz embedding two} and (\ref{uniform H21 bound ibu}), $$\norm{\nabla \theta_g}{L^{\infty}} \leq C_0 \norm{\theta_g}{W^{2,(2,1)}}\leq C_0K$$
with $C_0K$ independent of $g \in \mathcal{G}$. Thus, \eqref{L2 stab} holds with 
$
C_1=e^{CC_0KT}.
$
\end{proof}
Since the $W^{k,1}(\R^2)$ norm is not continuous with respect to the $L^2$-topology, we consider a functional that mollifies and localizes, and then show this functional is continuous in $L^2(\R^2)$. In particular, for a sufficient choice of mollification and localization parameters, the functional detects large $W^{k,1}(\R^2)$ norms.
\begin{lemma}
\label{functional}
Let $\varphi \in C_c^{\infty}(\R^2)$, with $\varphi \geq 0$ and $\int \varphi=1$, and for $h>0$ set $\varphi_h(x)=h^{-2}\varphi(x/h)$. For $R,h>0$, define $T_{R,h}: L^1_{\text{loc}}(\R^2) \to \R_{\geq 0}$ by 
$$
T_{R,h}(f) = \sum_{|\alpha| \leq k} \int_{B_R}|D^{\alpha}(\varphi_h \ast f)(x)|\,dx.
$$
For every $R,h>0$, the functional $T_{R,h}$ satisfies
\begin{equation}
\label{lemma 2 ibu 1}
T_{R,h}(f) \leq \norm{f}{W^{k,1}},
\qquad f\in W^{k,1}(\R^2).
\end{equation}
Moreover, there exists a constant $C_{R,h}>0$ such that
\begin{equation}
\label{lemma 2 ibu 2}
\left|T_{R,h}(f)-T_{R,h}(g)\right|
\leq C_{R,h}\norm{f-g}{L^2},
\qquad f,g\in L^2(\R^2).
\end{equation}
Finally, if $A>0$ and $f\in W^{k,1}(\R^2)$ satisfies
$\norm{f}{W^{k,1}}>4A$, then there exist $R,h>0$ such that
\begin{equation}
\label{lemma 2 ibu 3}
T_{R,h}(f)>3A.
\end{equation}
\end{lemma}
\begin{proof}
Let $f \in W^{k,1}(\R^2)$, since derivatives commute with convolution, from Young's convolution inequality it follows that 
$$
T_{R,h}(f) \leq \sum_{|\alpha| \leq k} \|\varphi_h \ast D^{\alpha}f\|_{L^{1}} \leq \|f\|_{W^{k,1}},
$$
where we used that $\|\varphi_h\|_{L^1}=1$. This gives (\ref{lemma 2 ibu 1}). To prove (\ref{lemma 2 ibu 2}), let $f,g \in L^2(\R^2)$. From the reverse triangle inequality, the Cauchy-Schwarz inequality, and Young's convolution inequality,
\begin{equation}
\begin{split}
|T_{R,h}(f)-T_{R,h}(g)| &\leq \sum_{|\alpha| \leq k} \int_{B_R}|D^{\alpha}(\varphi_h \ast (f-g))|\,dx\\
&\leq |B_R|^{1/2}\sum_{|\alpha| \leq k} \|D^\alpha\varphi_h\|_{L^1}\|f-g\|_{L^2},
\end{split}
\end{equation}
which gives \eqref{lemma 2 ibu 2} with 
$$
C_{R,h}=|B_R|^{1/2}\sum_{|\alpha|\leq k} \|D^{\alpha}\varphi_h\|_{L^1}.
$$
Finally, to prove (\ref{lemma 2 ibu 3}), let $f \in W^{k,1}(\R^2)$ satisfy $\|f\|_{W^{k,1}}  >4A$. Since $\varphi_h\ast f \to f$ in $W^{k,1}(\R^2)$ as $h \to 0^+$, we choose $h$ sufficiently small so that 
$$
\|\varphi_h \ast f - f\|_{W^{k,1}} < A/2.
$$
Then,
$$
\|\varphi_h \ast f\|_{W^{k,1}} \geq \|f\|_{W^{k,1}} - \|\varphi_h\ast f-f\|_{W^{k,1}} > 7A/2.
$$
With this $h$ fixed, we choose $R>0$ sufficiently large so that 
$$
\sum_{|\alpha| \leq k} \int_{\R^2 \setminus B_R} \left|D^{\alpha}(\varphi_h \ast f)\right|\,dx < A/2.
$$
Consequently
$$
T_{R,h}(f) = \|\varphi_h \ast f\|_{W^{k,1}}-\sum_{|\alpha| \leq k} \int_{\R^2 \setminus B_R} \left|D^{\alpha}(\varphi_h \ast f)\right|\,dx >3A.
$$

\end{proof}

\begin{lemma}
\label{big functional}
Let $A>0$, and let $T>0$ be as in Lemma \ref{H21 LWP and L2 stab}. Suppose $f \in \mathcal{G} \cap C_c^{\infty}(\R^2)$ and its corresponding solution $\theta_f$ to \eqref{SQG} is smooth through $[0,T]$. Also assume that for some time $t \in (0,T)$
$$
\norm{\theta_f(t)}{W^{k,1}} > 4A.
$$
Then, there exist $R,h,\rho>0$ such that whenever $g \in \mathcal{G}$ and 
\begin{equation}
\label{lemma 3 ibu 1}
\norm{g-f}{W^{k,1}}<\rho,
\end{equation} the corresponding unique weak solution $\theta_g$ satisfies
\begin{equation}
\label{lemma 3 ibu 2}
T_{R,h}(\theta_g(t))>2A.
\end{equation}
\end{lemma}
\begin{proof}
Since $\|\theta_f(t)\|_{W^{k,1}} >4A,$ from Lemma \ref{functional} there exist $R,h>0$ such that 
$$
T_{R,h}(\theta_f(t))>3A.
$$
Now, given $g \in \mathcal{G}$, from Lemma \ref{H21 LWP and L2 stab}, its unique weak solution exists on $[0,T]$. Furthermore,  from (\ref{lemma 2 ibu 2}), Lemma \ref{H21 LWP and L2 stab}, and Lemma \ref{functional} 
\begin{equation}
\begin{split}
|T_{R,h}(\theta_g(t))-T_{R,h}(\theta_f(t))| &\leq C_{R,h}\exp\left(\int_0^T \|\nabla \theta_g(\tau)\|_{L^{\infty}}\,d\tau\right)\|f-g\|_{L^2}\\
&\leq CC_{R,h}C_1\|f-g\|_{W^{k,1}},
\end{split}
\end{equation}
where $C_1$ is as in \eqref{L2 stab}.
Hence, we choose 
$$
\rho =  \frac{A}{CC_{R,h}C_1}.
$$
Then, whenever $g \in \mathcal{G}$ satisfies $\|f-g\|_{W^{k,1}}<\rho$
$$
|T_{R,h}(\theta_g(t))-T_{R,h}(\theta_f(t))| \leq A.
$$
It immediately follows that 
$$
T_{R,h}(\theta_g(t)) \geq T_{R,h}(\theta_f(t)) - |T_{R,h}(\theta_f(t))-T_{R,h}(\theta_g(t))| >3A-A=2A.
$$
\end{proof}

\begin{proposition}
\label{initial sequence}
There exists a sequence of initial data $\theta_0^{(n)} \in \mathcal{G} \cap C_c^{\infty}$, a sequence $t_n \downarrow 0,$ and sequences of parameters $R_n,h_n,\rho_n>0$ such that 
\begin{equation}
\label{prop 1 ibu 1}
\norm{\theta^{(n)}(t_n)}{W^{k,1}}>4n,
\end{equation}
and whenever $g \in \mathcal{G}$ with $\norm{g-\theta_0^{(n)}}{W^{k,1}}<\rho_n$, its corresponding solution $\theta_g$ to \eqref{SQG} satisfies
\begin{equation}
\label{prop 1 ibu 2}
T_{R_n,h_n}(\theta_g(t_n))>2n.
\end{equation}
Moreover, if $\supp \overline{\theta_0} \subseteq B_L(0)$ for some $L>0$, then for all $n\in\mathbb{N}$,  $\theta_0^{(n)}$ is supported in $B_{L}(0)$ and 
\begin{equation}
\label{prop 1 ibu 3}
\norm{\theta_0^{(n)}-\theta_0^{(n-1)}}{W^{k,1}} \leq 2^{-(n+2)} \min\{r,\rho_1,\cdots,\rho_{n-1}\},
\end{equation}
where $\theta_0^{(0)}=\overline{\theta}_0$ and the minimum is given by $r$ when $n=1$.
\end{proposition}
\begin{proof}
Let $T>0$ be as in Lemma \ref{H21 LWP and L2 stab}, and set $\theta_0^{(0)} = \overline{\theta}_0$ and $t_0=T$.  Assume $L>0$ is such that $\supp \overline{\theta}_0 \subseteq B_L(0)$. We proceed by induction. We first construct $\theta_0^{(1)}$. Choose $t_1$ satisfying
$$
0<t_1<\min\{2^{-1},t_0/2\}.
$$
Apply Theorem \ref{norm inflation key} with background solution $\overline{\theta},$ time $\tau=t_1$, and parameters
$$
\varepsilon=2^{-3}r, \qquad M=4+\|\overline{\theta}(t_1)\|_{W^{k,1}},
$$
to obtain $\theta_0^{(1)} \in C_c^\infty(\R^2)$ with $\supp\theta_0^{(1)} \subseteq B_L(0)$ such that 
$$
\|\theta_0^{(1)}-\theta_0^{(0)}\|_{W^{k,1}}<2^{-3}r
$$
and 
$$
\|\theta^{(1)}(t_1)-\overline{\theta}(t_1)\|_{W^{k,1}} >4 +\|\overline{\theta}(t_1)\|_{W^{k,1}}.
$$
The reverse triangle inequality gives $\|\theta^{(1)}(t_1)\|_{W^{k,1}}>4$. Furthermore, since $\|\theta_0^{(1)}-\overline{\theta}_0\|_{W^{k,1}}<r$, we have $\theta_0^{(1)} \in \mathcal{G}$. We then apply Lemma \ref{big functional} with $f=\theta_0^{(1)},$ $A=1$, and $t=t_1$ to obtain $R_1,h_1,\rho_1>0$ such that every $g \in \mathcal{G}$ satisfying 
$$
\|g-\theta_0^{(1)}\|_{W^{k,1}}<\rho_1,
$$
also satisfies
$$
T_{R_1,h_1}(\theta_g(t_1))>2,
$$
establishing the base case.

Now, suppose that for some $ n \geq 2$ the construction has been completed to obtain $\theta_0^{(j)}, t_j,R_j,h_j,$ and $\rho_j$ for $1 \leq j \leq n-1$. Note that since $\theta^{(0)}_0 = \overline{\theta}_0$,
 $$
 \|\theta_0^{(n-1)} - \overline{\theta}_0\|_{W^{k,1}} \leq \sum_{j=1}^{n-1}\|\theta_0^{(j)} - \theta_0^{(j-1)}\|_{W^{k,1}} \leq \frac{r}{4}\sum_{j=1}^{n-1} 2^{-j} < \frac{r}{4}.
 $$
 Hence, from our choice of $r$, we see that $\theta_0^{(n-1)}\in \mathcal{G}\cap C_c^{\infty}(\R^2)$ is non-constant. Furthermore, from Lemma \ref{H21 LWP and L2 stab}, we see that if $\theta^{(n-1)}$ is the corresponding unique weak solution, then $\|\nabla u^{(n-1)}(t)\|_{L^{\infty}}$ and $\|\nabla \theta^{(n-1)}(t)\|_{L^{\infty}}$ are uniformly bounded on $[0,T]$. In particular, as this quantity controls all higher-order Sobolev norms, each $\theta^{(n-1)}$ remains smooth on a uniform time interval $[0,T]$. We choose $t_n$ satisfying 
 $$
 0<t_n<\min\left\{ 2^{-n}, \frac{1}{2}t_{n-1}\right\}.
 $$
 Then $t_n \downarrow 0$. We again apply Theorem \ref{norm inflation key} with 
 $$
 \tau=t_n, \quad \varepsilon = 2^{-(n+2)} \min\left\{ r,\rho_1,\cdots,\rho_{n-1} \right\},\quad  M=4n + \|\theta^{(n-1)}(t_n)\|_{W^{k,1}}
 $$
 to obtain $\theta_0^{(n)} \in C_c^{\infty}(\R^2)$ with $\supp \theta_0^{(n)} \subseteq B_L(0)$,
 $$
 \|\theta_0^{(n)}-\theta_0^{(n-1)}\|_{W^{k,1}}<\varepsilon,
 $$
 and 
 $$
 \|\theta^{(n)}(t_n)-\theta^{(n-1)}(t_n)\|_{W^{k,1}} >4n+\|\theta^{(n-1)}(t_n)\|_{W^{k,1}}.
 $$
 Note from the reverse triangle inequality
 $$
 \|\theta^{(n)}(t_n)\|_{W^{k,1}} >4n.
 $$
 Furthermore, 
 $$
 \|\theta_0^{(n)}-\overline{\theta_0}\|_{W^{k,1}} \leq \frac{r}{4}\sum_{j=1}^n 2^{-j} <\frac{r}{4}<r,
 $$
 so that $\theta_0^{(n)} \in \mathcal{G}$. We apply Lemma \ref{big functional} with $f=\theta_0^{(n)}$, $A=n$, $t=t_n$ to supply $R_n,h_n,\rho_n>0$ such that every $g \in \mathcal{G}$ satisfying 
 $$
 \|g-\theta_0^{(n)}\|_{W^{k,1}} <\rho_n,
 $$
 also satisfies 
 $$
 T_{R_n,h_n}(\theta_g(t_n)) >2n.
 $$
\end{proof}
\subsection{Proof of Theorem \ref{BLOW UP}: Instantaneous Blow-up}
We now pass to the limit in the sequence constructed in Proposition \ref{initial sequence}. The summable increments ensure convergence to an initial datum $\theta_0$ in $W^{k,1}(\R^2)$, while the increments' dependence on the radii $\rho_n$ guarantees that 
$$
\|\theta_0-\theta_0^{(n)}\|_{W^{k,1}}<\rho_n \quad \text{for every }n\in\mathbb{N}.
$$
 Consequently, the solution $\theta$ generated from $\theta_0$ inherits the lower bound $T_{R_n,h_n}(\theta(t_n))>2n$. The functionals $T_{R_n, h_n}$ are bounded above by the $W^{k,1}$ norm, yielding the desired instantaneous blow-up. We include the details of this argument in the proof of Theorem \ref{initial data} below, from which Theorem \ref{BLOW UP} follows.
\begin{theorem}
\label{initial data}
Let $\theta_0^{(n)}$, $\rho_n$, and $L$ be as in Proposition \ref{initial sequence}. Then there exists some $\theta_0 \in W^{k,1}(\R^2)$ such that $\theta_0^{(n)} \to \theta_0$ in $W^{k,1}(\R^2)$. Moreover,
\begin{equation}
\label{prop 2 ibu 1}
\norm{\overline{\theta}_0-\theta_0}{W^{k,1}} <\eps, \quad \supp \theta_0 \subset  \overline{B_L(0)},
\end{equation}
and, for every $n \in \mathbb{N}$,
\begin{equation}
\label{prop 2 ibu 2}
\norm{\theta_0 - \theta_0^{(n)}}{W^{k,1}} <\rho_n.
\end{equation}
Consequently, the unique solution $\theta$ generated by $\theta_0$ as given by Theorem \ref{LWP H21} on $[0,T]$ satisfies
\begin{equation}
\label{prop 2 ibu 3}
\esssup_{0<t<T'}\|\theta(t)\|_{W^{k,1}}=\infty, \quad \text{for every } 0<T'\leq T.
\end{equation}
\end{theorem}
\begin{proof}
It follows from (\ref{prop 1 ibu 3}) that for all $j\geq 1$,
$$
\|\theta_0^{(j)}-\theta_0^{(j-1)}\|_{W^{k,1}} \leq r2^{-(j+2)}.
$$
Hence, for $m>n\geq 0$
\begin{equation}
\label{summable so cauchy}
\|\theta_0^{(m)}-\theta_0^{(n)}\|_{W^{k,1}} \leq \sum_{j=n+1}^m \|\theta_0^{(j)}-\theta_0^{(j-1)}\|_{W^{k,1}} \leq r \sum_{j=n+1}^\infty 2^{-(j+2)} =r2^{-(n+2)}.
\end{equation}
Thus $\theta_0^{(n)}$ is Cauchy in $W^{k,1}(\R^2)$. By completeness, there exists some $\theta_0 \in W^{k,1}(\R^2)$ such that $\theta_0^{(n)} \to \theta_0$ in $W^{k,1}(\R^2)$. Furthermore, taking $n=0$ in \eqref{summable so cauchy} and sending $m \to \infty$ yields 
$$
\|\theta_0 - \overline{\theta}_0\|_{W^{k,1}} \leq \frac{r}{4} <\eps,
$$
so that $\theta_0 \in \mathcal{G}$, and the first statement in (\ref{prop 2 ibu 1}) holds.

Arguing in a similar way, since for $j \geq n+1$ 
$$
\|\theta_0^{(j)} - \theta_0^{(j-1)}\|_{W^{k,1}} \leq 2^{-(j+2)}\rho_n,
$$
for $m>n$ one has 
\begin{equation}
\label{sending m to infinity again}
\|\theta_0^{(m)} - \theta_0^{(n)}\|_{W^{k,1}} \leq \rho_n\sum_{j=n+1}^m 2^{-(j+2)} \leq 2^{-(n+2)}\rho_n <\rho_n.
\end{equation}
Sending $m \to \infty$ in \eqref{sending m to infinity again} yields
$$
\|\theta_0-\theta_0^{(n)}\|_{W^{k,1}} <\rho_n,
$$
giving (\ref{prop 2 ibu 2}).

Finally, to see that $\supp \theta_0 \subseteq \overline{B_L(0)}$, since $k \geq 3$, the embedding $W^{k,1}(\R^2) \hookrightarrow BUC(\R^2)$ implies that $\theta_0^{(n)} \to \theta_0$ uniformly. As each $\theta_0^{(n)}$ vanishes outside $B_L(0),$ so does $\theta_0$. Consequently, 
$$
\supp \theta_0 \subseteq \overline{B_L(0)}.
$$
To complete the proof of Theorem \ref{initial data}, it remains to establish (\ref{prop 2 ibu 3}). Since $\|\theta_0-\theta_0^{(n)}\|_{W^{k,1}} <\rho_n$ for every $n \in \mathbb{N}$, from Proposition \ref{initial sequence}, it follows that $T_{R_n,h_n}(\theta(t_n)) >2n$. 
For a fixed $n$ and $0<t<T$, Lemma \ref{functional} gives 
$$
|T_{R_n,h_n}(\theta(t))-T_{R_n,h_n}(\theta(t_n))| \leq C_{R_n,h_n}\|\theta(t)-\theta(t_n)\|_{L^2}.
$$
Since $\theta \in C([0,T];W^{1,(2,1)}(\R^2))$, the above right hand side goes to $0$ as $t \to t_n$ (as $t_n$ is fixed). Hence, there exists some $\delta_n>0$ such that whenever $|t-t_n|<\delta_n$
$$
\left|T_{R_n,h_n}(\theta(t))-T_{R_n,h_n}(\theta(t_n))\right| <n.
$$ 
If $\theta(t) \in W^{k,1}(\R^2)$, then, from the reverse triangle inequality
$$
\|\theta(t)\|_{W^{k,1}} \geq T_{R_n,h_n}(\theta(t)) > T_{R_{n},h_n}(\theta(t_n))-n>n,
$$
for all $t<T$ satisfying $|t-t_n|<\delta_n$. Otherwise, $\|\theta(t)\|_{W^{k,1}} = \infty$ and the same lower bound holds. Now, fix $T' \in (0,T]$ and $M>0$. Since $t_n \downarrow 0$, we can choose $n \in \mathbb{N}$ satisfying $$n>M\quad\text{and} \quad0<t_n<T'.$$
Shrinking $\delta_n$ if necessary, we have that for each $t \in (t_n-\delta_n,t_n+\delta_n) \subseteq   (0,T')$,
$$
\|\theta(t)\|_{W^{k,1}}>M.
$$
Since $M>0$ is arbitrary and this lower bound is detected on a set of positive measure, \eqref{prop 2 ibu 3} follows. Uniqueness of $\theta(t)$ follows from Theorem \ref{LWP H21}. 
\end{proof}

\begin{appendix}
\section{Energy Estimates} 
\label{appendix energy}
In this section, we collect and prove the energy estimates used to establish local norm-inflation. Throughout the section, $k\geq 3$ is fixed, and $\overline{\theta}$, $\overline{u}$, $\tau$, and $\Psi$ are defined as in Section \ref{norm inflation section}. In addition, as in Section \ref{norm inflation section}, for $r>0$ we set
\begin{equation*}
       \begin{split}
           U_r(t,y)  = \frac{\overline{u}(t,ry)}{r},\quad G_r (t,y)= \nabla\overline{\theta}(t,ry),
       \end{split}
   \end{equation*}
   and 
   \begin{equation*}
       E = \partial_t\Psi+U_r\cdot\nabla \Psi +G_r\cdot R^\perp\Psi.
   \end{equation*}
\begin{lemma}
\label{energy estimate 1}
Suppose that $a>0$, $0 < r \leq 1$, and $v$ is a smooth solution to 
\begin{equation}
\label{energy estimate 1 pde}
\begin{cases}
\partial_tv+R^{\perp}v\cdot G_r(t,y) + U_r(t,y) \cdot \nabla v + ar^{k-3}R^{\perp}v\cdot \nabla v=0,\\
v|_{t=0} = v_0,
\end{cases}
\end{equation}
with $\nabla\cdot(U_r) = 0$, and assume there exists $M>0$ independent of $r$ such that
\begin{equation}\label{overlinebounds}
\sup_{0 \leq t \leq \tau} \left(\sum_{1 \leq |\beta| \leq k} \norm{ D^{\beta}U_r(t)}{L^{\infty}} + \norm{G_r(t)}{W^{k,\infty}}\right) \leq M.
\end{equation}
Then there exists a $C>0$ depending on $M$ and $k$, but independent of $a,r$ and $v_0$, such that 
\begin{equation}
\label{energy estimate 1 conclusion}
 \|v(t)\|_{H^k} \leq \frac{e^{Ct} \left\| v_0\right\|_{H^k}}{1-Ca\left\|v_0\right\|_{H^k}\int_0^t e^{Cs}\,ds},
\end{equation}
whenever $t\in [0,\tau]$ is sufficiently small so that $1-Ca\left\|v_0\right\|_{H^k}\int_0^t e^{Cs}\,ds>0$.
\end{lemma}
\begin{remark}
To see why $M$ in \eqref{overlinebounds}  is independent of $r$ in our application, from the chain rule we have 
$$
D_y^{\beta}U_r(t,y)=r^{|\beta|-1}(D^{^\beta}\overline{u})(t,ry), \quad 1\leq |\beta| \leq k,$$
and 
$$
D^{\beta}_yG_r(t,y)= r^{|\beta|}(D^\beta \nabla \overline{\theta})(t,ry), \quad 0 \leq |\beta|\leq k.
$$
Since $0<r\leq 1$ and the background solution is smooth, the above quantities are bounded in $L^{\infty}(\R^2)$, uniformly in $t \in [0,\tau]$ and $r \in (0,1]$. Thus $M$ may be chosen to depend only on $k$ and the background solution on $[0,\tau]$.
\end{remark}
\begin{proof}
Let $\alpha$ be a multi-index with $|\alpha| \leq k$, and fix $t\in [0,\tau]$. Since $U_r$ and $R^{\perp}v$ are divergence free, 
\begin{equation}
\label{divergence free energy estimate 1}
\int_{\R^2} U_r(t,x) \cdot \nabla  D^{\alpha}v D^{\alpha}v\,dx=0 \quad \text{and } \int_{\R^2} R^{\perp}v\cdot \nabla  D^{\alpha}v D^{\alpha}v\,dx=0.
\end{equation}
Applying $ D^{\alpha}$ to \eqref{energy estimate 1 pde}, taking an inner-product with $ D^{\alpha}v$, and using \eqref{divergence free energy estimate 1} we see that 
\begin{equation}
    \label{ee1 eq1}
    \begin{split}
    &\qquad \frac{1}{2}\frac{d}{dt}\norm{v(t)}{H^k}^2 = - \sum_{|\alpha|\leq k}\int_{\R^2}\left[ D^{\alpha},U_r(t,x)\cdot \nabla\right]v D^{\alpha}v\,dx \\
     &- \sum_{|\alpha|\leq k}\int_{\R^2} D^{\alpha}(G_r(t,x)\cdot R^{\perp}v) D^{\alpha}v\,dx  -\sum_{|\alpha|\leq k}ar^{k-3}\int_{\R^2} [ D^{\alpha},R^{\perp}v\cdot\nabla]v D^{\alpha}v\,dx.
    \end{split}
\end{equation}
We estimate each term on the right hand side of (\ref{ee1 eq1}). For the first term, we apply the Leibniz rule and (\ref{overlinebounds}), which gives 
\begin{equation}
\label{ee1 eq2}
\begin{split}
&\qquad\quad \left |\sum_{|\alpha|\leq k}\int_{\R^2}\left[ D^{\alpha},U_r(t,x)\cdot \nabla\right]v D^{\alpha}v\,dx\right |\\
&\leq C\norm{v}{H^k}\sum_{|\alpha|\leq k}\sum_{\{ |\beta|>0: \beta \leq \alpha\}} C_{\alpha,
\beta} \norm{D^{\beta}U_r(t) \cdot \nabla  D^{\alpha-\beta}v}{L^2} \leq CM\norm{v}{H^k}^2.
\end{split}
\end{equation}
To estimate the second term on the right hand side of (\ref{ee1 eq1}), we apply the Cauchy-Schwarz inequality, the Leibniz rule, boundedness of Riesz operators on $H^k(\R^2)$, and the bound on $G_r$ in (\ref{overlinebounds}). We conclude that  
\begin{equation}
\label{ee1 eq3}
 \left |\sum_{|\alpha|\leq k} \int_{\R^2} D^{\alpha}(G_r(t,x)\cdot R^{\perp}v) D^{\alpha}v\,dx\right | \leq CM\norm{v}{H^k}^2.
 \end{equation}
Finally, for the third term on the right hand side of (\ref{ee1 eq1}), we first note that, by the standard commutator estimate 
\begin{equation}\label{commutatorest}
\| [D^{\alpha},f\cdot\nabla] g \|_{L^2} \leq C(\|\nabla f \|_{L^{\infty}}\| g \|_{H^k} + \|\nabla g\|_{L^{\infty}}\| f \|_{H^k}),
\end{equation}
for $f,g\in H^k(\R^2)$, combined with the Sobolev embedding theorem and boundedness of Riesz operators on $H^k(\R^2)$,
\begin{equation}
\label{nonlinear commutator estimate}
\begin{split}
&\sum_{|\alpha|\leq k}
\left\|[D^\alpha,R^\perp v\cdot\nabla]v\right\|_{L^2}\leq C\left(
\|\nabla R^\perp v\|_{L^\infty}\|v\|_{H^k}
+\|\nabla v\|_{L^\infty}\|R^\perp v\|_{H^k}
\right)\leq C\|v\|_{H^k}^2.
\end{split}
\end{equation}
The Cauchy-Schwarz inequality, the bound $ar^{k-3} \leq a$ for $k \geq 3$ and $0 < r \leq 1$, and (\ref{nonlinear commutator estimate}) allow us to write
\begin{equation}
\label{ee1 eq4}
\begin{split}
&\left |ar^{k-3}\sum_{|\alpha|\leq k}\int_{\R^2} [ D^{\alpha},R^{\perp}v\cdot\nabla]v D^{\alpha}v\,dx \right | \leq a \norm{v}{H^k} \sum_{|\alpha| \leq k}\norm{[ D^{\alpha},R^{\perp}v\cdot\nabla]v}{L^2}\\
&\qquad \leq Ca\norm{v}{H^k}^2\left(\norm{\nabla R^{\perp}v}{L^{\infty}} + \norm{\nabla v}{L^{\infty}}\right) \leq Ca\norm{v}{H^k}^3.
\end{split}
\end{equation}
Substituting \eqref{ee1 eq2}, \eqref{ee1 eq3}, and \eqref{ee1 eq4} into \eqref{ee1 eq1}, it follows that 
\begin{equation}
\label{ee1 eq5}
\frac{d}{dt}\norm{v(t)}{H^k} \leq C\norm{v(t)}{H^k} +Ca\norm{v(t)}{H^k}^2,
\end{equation}
where $C$ now depends on $M$. Solving the differential inequality, we find that
$$
\norm{v(t)}{H^k} \leq \frac{e^{Ct}\norm{v_0}{H^k}}{1-Ca\norm{v_0}{H^k}\int_0^t e^{Cs}\,ds}
$$
for all $t\in [0,\tau]$ for which the denominator remains positive.
\end{proof}
\begin{lemma}
\label{energy estimate 2}
Let $a>0$, $0 < r \leq 1$, and assume $w$ is a smooth solution to
\begin{equation}
\label{energy estimate 2 pde}
\begin{split}
&\partial_tw+(U_r+ar^{k-3}R^{\perp}v)\cdot\nabla w + G_r\cdot R^{\perp}w \\
&\qquad + ar^{k-3}R^{\perp}w\cdot \nabla \Psi =-E - ar^{k-3}R^{\perp}\Psi\cdot \nabla \Psi
\end{split}
\end{equation}
with $\nabla\cdot(U_r)=0$, and with $U_r$ and $G_r$ satisfying (\ref{overlinebounds}). Then there exist $C>0$ depending on $M$ and $k$, but independent of $r$, $a$, $v$, and $\Psi$ such that
\begin{equation}
\label{energy estimate 2 conclusion}
\begin{split}
\frac{d}{dt}\norm{w(t)}{H^k} &\leq C\left(1+ar^{k-3}\norm{v(t)}{H^k}+ar^{k-3}\norm{\Psi(t)}{H^{k+1}}\right)\norm{w(t)}{H^k} + \norm{E(t)}{H^k}\\ &\quad+ Car^{k-3}\norm{\Psi(t)}{H^k}\norm{\Psi(t)}{H^{k+1}}.
\end{split}
\end{equation}
\end{lemma}
\begin{proof}
For each multi-index $\alpha $ with $|\alpha|\leq k$, we apply $ D^{\alpha}$ to \eqref{energy estimate 2 pde}, take an $L^2$ inner-product with $ D^{\alpha}w$, and use that $U_r + ar^{k-3}R^{\perp}v$ is divergence free to write
\begin{equation}
\label{ee2 eq 1}
\begin{split}
\frac{1}{2}\frac{d}{dt}\norm{w}{H^k}^2 =& -\sum_{|\alpha| \leq k} \langle [ D^{\alpha},(U_r+ar^{k-3}R^{\perp}v)\cdot \nabla]w, D^{\alpha}w \rangle-\sum_{|\alpha|\leq k} \langle  D^{\alpha} (G_r\cdot R^{\perp}w),  D^{\alpha}w\rangle\\
&\hspace{-2cm}-ar^{k-3}\sum_{|\alpha|\leq k} \langle  D^{\alpha}(R^{\perp}w\cdot \nabla \Psi), D^{\alpha}w\rangle-\sum_{|\alpha| \leq k}\langle  D^{\alpha}E, D^{\alpha}w\rangle 
-ar^{k-3}\sum_{|\alpha| \leq k} \langle  D^{\alpha}\left(R^{\perp}\Psi\cdot \nabla \Psi\right), D^{\alpha}w\rangle.
\end{split}
\end{equation}
We estimate each term on the right-hand side, beginning with the first term. We estimate the commutator involving $U_r$ by the Leibniz rule and \eqref{overlinebounds}, as we did in \eqref{ee1 eq2}. For the commutator involving $R^{\perp}v$, we apply (\ref{commutatorest}) combined with the Sobolev embedding theorem, and the boundedness of Riesz operators on $H^k(\R^2)$ to obtain
\begin{equation}
\label{ee2 eq 2}
\begin{split}
&\qquad \left |\sum_{|\alpha| \leq k} \langle [ D^{\alpha},(U_r+ar^{k-3}R^{\perp}v)\cdot \nabla]w, D^{\alpha}w \rangle\right | \\
&\leq \norm{w}{H^k}\sum_{|\alpha|\leq k } \left(\norm{[ D^{\alpha},U_r\cdot \nabla]w} {L^2} + ar^{k-3}\norm{[ D^{\alpha},R^{\perp}v\cdot \nabla]w}{L^2}\right)\\
& \leq C\norm{w}{H^k}\left(\norm{w}{H^k} + ar^{k-3}\norm{\nabla R^{\perp}v}{L^{\infty}}\norm{w}{H^k} + ar^{k-3}\norm{R^{\perp}v}{H^k}\norm{\nabla w}{L^{\infty}}\right)\\ 
&\leq C\norm{w}{H^k}^2(1+ar^{k-3}\norm{v}{H^k}),
\end{split}
\end{equation}
where $C$ now depends on $M$. For the second term on the right hand side of (\ref{ee2 eq 1}), we again apply the Cauchy-Schwarz inequality, Leibniz rule, boundedness of Riesz operators on $H^k(\R^2)$, and (\ref{overlinebounds}) to write 
\begin{equation}
\label{ee2 eq 3}
\left |\sum_{|\alpha|\leq k} \langle  D^{\alpha} (G_r\cdot R^{\perp}w), D^{\alpha}w\rangle \right | \leq C\norm{w}{H^k}^2.  
\end{equation}
We can apply similar methods to those in (\ref{ee2 eq 3}), as well as the observation that $H^k(\R^2)$ is a Banach algebra for $k \geq 3$,  to bound the remaining terms on the right hand side of (\ref{ee2 eq 1}).  We have  
\begin{equation}
\label{ee2 eq 4}
\begin{split}
&\abs{ar^{k-3}\sum_{|\alpha|\leq k} \langle  D^{\alpha}(R^{\perp}w\cdot \nabla \Psi), D^{\alpha}w\rangle } \leq Car^{k-3} \norm{w}{H^k}^2\norm{\Psi}{H^{k+1}},\\
&\abs{\sum_{|\alpha| \leq k}\langle  D^{\alpha}E, D^{\alpha}w\rangle } \leq \norm{E}{H^k}\norm{w}{H^k}, \\
&\abs{ar^{k-3}\sum_{|\alpha| \leq k} \langle  D^{\alpha}\left(R^{\perp}\Psi\cdot \nabla \Psi\right), D^{\alpha}w\rangle  } \leq Car^{k-3}\norm{\Psi}{H^k}\norm{\Psi}{H^{k+1}}\norm{w}{H^k}.
\end{split}
\end{equation}
Substituting \eqref{ee2 eq 2}, \eqref{ee2 eq 3}, and \eqref{ee2 eq 4} into \eqref{ee2 eq 1} yields \eqref{energy estimate 2 conclusion}.
\end{proof}

\section{Local Well-Posedness in $W^{2,(2,1)}(\R^2)$}
\label{LWP}
\subsection{Sobolev-Lorentz Spaces}
We begin by recalling the definition of Lorentz spaces. Given $f:\R^d \to \R$, we first define the \emph{distribution function}
$$
d_f(\alpha)=\left|\left\{ x \in \R^d : |f(x)|>\alpha\right\}\right|,
$$
and the \emph{decreasing rearrangement} of $f$, denoted $f^*:[0,\infty) \to [0,\infty]$ as 
$$
f^*(t) = \inf\left\{s>0: d_f(s) \leq t \right\}.
$$
\begin{definition} For $p\in[1,\infty)$ and $q \in [1,\infty]$, we define the \emph{Lorentz space} $L^{(p,q)}(\R^d)$ to be the space of measurable functions satisfying $\norm{f}{L^{(p,q)}}<\infty$, where
\begin{equation}
\label{Lorentz space definition}
\norm{f}{L^{(p,q)}} =
\begin{cases}
\left(\int_0^\infty \left(t^{\frac{1}{p}}f^*(t)\right)^q\frac{dt}{t}\right)^{\frac{1}{q}} \quad &\text{if } q<\infty,\\
\sup_{t>0} t^\frac{1}{p} f^*(t) & \text{if } q=\infty.
\end{cases}
\end{equation}
For $k \geq 0$ an integer, we also define the \emph{Sobolev-Lorentz space} $W^{k,(p,q)}(\R^d)$ as 
$$
W^{k,(p,q)}(\R^d) = \left\{ f \in L^{(p,q)}(\R^d) : \norm{f}{W^{k,(p,q)}}<\infty\right\},
$$
where 
$$
\norm{f}{W^{k,(p,q)}} = \sum_{|\alpha| \leq k} \norm{D^{\alpha}f}{L^{(p,q)}}.
$$
\end{definition}
Under the above norms, both $W^{k,(p,q)}(\R^d)$ and $L^{(p,q)}(\R^d)$ are complete, $L^{(p,q)}(\R^d) \hookrightarrow L^{(p,r)}$ for $q<r$, and $L^{(p,p)}(\R^d)=L^p(\R^d)$ for all $p \in [1,\infty)$. 

Lorentz spaces enjoy many of the same properties as $L^p$ spaces, including analogues of H\"{o}lder's and Young's inequalities, which we recall below. 
\begin{proposition}[\cite{Hunt66}]\label{Lorentz Holder} Let $1 < p,p_1,p_2<\infty, 1\leq q,q_1,q_2 \leq \infty$ satisfy $\frac{1}{p}=\frac{1}{p_1}+\frac{1}{p_2}$ and $\frac{1}{q} = \frac{1}{q_1}+\frac{1}{q_2}$. Then there exists a constant $C$ (depending only on $p_1,p_2,q_1$ and $q_2$) such that for any $f \in L^{(p_1,q_1)}(\R^d)$ and $g \in L^{(p_2,q_2)}(\R^d)$ 
$$
\norm{fg}{L^{(p,q)}} \leq C\norm{f}{L^{(p_1,q_1)}}\norm{g}{L^{(p_2,q_2)}}.
$$
Moreover, if $g\in L^{(p,q)}(\R^d)$ and $f\in L^{\infty}(\R^d)$, then 
$$ \|fg \|_{L^{(p,q)}} \leq \|f\|_{L^{\infty}} \|g \|_{L^{(p,q)}}.$$
\end{proposition}
Next, we recall a convolution inequality on Lorentz spaces. This inequality first appeared in \cite{ONeil63}; however, that result does not cover the more general case where one of the functions in the convolution belongs to $L^1(\R^d)$. We state this case in the following proposition. The proof follows real interpolation methods and an application of Minkowski's integral inequality. For details, we refer the reader to Lemma 2.1 in \cite{GP22}. 
\begin{proposition}\label{convolution on lorentz} Let $1<p<\infty$ and $1 \leq q \leq p$. Suppose $f \in L^1(\R^d)$ is nonnegative and $g \in L^{(p,q)}(\R^d)$.  Then
$$
\|f\ast g\|_{L^{(p,q)}} \leq \|f\|_{L^1}\|g\|_{L^{(p,q)}}.
$$
\end{proposition}

We also have the following boundedness property of Calder\'on-Zygmund operators on Lorentz spaces.
\begin{proposition}[\cite{Yoshinaga69}]\label{CZO lorentz}
Let $T$ be a Calder\'{o}n-Zygmund operator, $1 < p<\infty$, and $1 \leq q \leq \infty$. Then there exists a constant $C$, depending only on $T,p,q$, such that $$\norm{Tf}{L^{(p,q)}} \leq C\norm{f}{L^{(p,q)}}.$$
\end{proposition}
One advantage of Lorentz spaces over the classical Lebesgue spaces is that they offer sharper Sobolev embeddings. We refer the reader to Lemma 31.1 in \cite{Tartar07} for a modern reference.
\begin{theorem}[Lemma 31.1, \cite{Tartar07}]\label{lorentz embedding 1}  Let $d \geq 2$, then $W^{1,1}(\R^d)$ is continuously embedded into the Lorentz space $L^{\left(\frac{d}{d-1},1\right)}(\R^d)$.
\end{theorem}
\noindent From Theorem \ref{lorentz embedding 1}, we see that $W^{1,1}(\R^2) \hookrightarrow L^{(2,1)}(\R^2)$, which is indeed sharper than the standard endpoint Sobolev embedding $W^{1,1}(\R^2) \hookrightarrow L^2(\R^2)$.

In Theorem \ref{lorentz embedding two} below, we recall the endpoint embedding $$W^{1,(d,1)}(\R^d) \hookrightarrow L^{\infty}(\R^d).$$ (In fact, one can replace $L^{\infty}(\R^d)$ above with $C_0(\R^d)$, the set of continuous functions vanishing at infinity). The proof of the theorem is not new, but we include it for completeness.
\begin{theorem}[\cite{Stein81}]
\label{lorentz embedding two}
Let $d \geq 1$. Then there exists a constant $C$, depending only on $d$ such that 
$$
\norm{f}{L^{\infty}} \leq C_d\norm{\nabla f}{L^{(d,1)}}, 
$$
for all $f \in W^{1,(d,1)}(\R^d)$.
\end{theorem}
\begin{proof} 
Let $f \in C_c^{\infty}(\R^d)$. We consider two cases. If $d=1$, $L^{(1,1)}(\R)=L^1(\R)$, and an application of the fundamental theorem of calculus yields $\|f\|_{L^\infty} \leq \|f'\|_{L^1}$. If $d \geq 2$, in the sense of distributions, one has that 
$$
\text{div}_y\left(\frac{y-x}{|y-x|^d}\right)=\sigma_{d-1}\delta_x(y),
$$
where $\sigma_{d-1}$ is the surface area of the unit sphere in $\R^d$, and $\delta_x$ is the Dirac delta distribution centered at $x$. Then, integrating by parts and using duality (see \cite{ONeil63}) and translation invariance of Lorentz norms, one has
\begin{equation*}
\begin{split}
|f(x)|&=\frac{1}{\sigma_{d-1}}\left|\int_{\R^d}  f(y)\text{div}_y\left(\frac{y-x}{|y-x|^d}\right)\,dy\right|\\ &\leq \frac{1}{\sigma_{d-1}} \norm{\nabla f}{L^{(d,1)}}\norm{\frac{1}{|\cdot|^{d-1}}}{L^{(d/(d-1),\infty)}}\\
&\leq C_d\norm{\nabla f}{L^{(d,1)}}.
\end{split}
\end{equation*}
The general case follows from the density of $C^\infty_{c}(\R^d)$ in $W^{1,(d,1)}(\R^d)$.
\end{proof}
\subsection{Proof of Local Well-Posedness}
Now that we have introduced some of the elementary properties of Lorentz spaces, we are equipped to establish local well-posedness of ($SQG$) in $W^{2,(2,1)}(\R^2)$. First, we define what we mean by a weak solution to \eqref{SQG}.
\begin{definition}
\label{weak solution definition}
We say that $\theta \in L^2([0,T]\times\R^2)$ is a weak solution to \eqref{SQG} with initial data $\theta_0$ if, setting $u=R^{\perp}\theta$, the pair ($u,\theta$) satisfies 
$$
\int_0^T \int_{\R^2} (\theta \partial_t\phi + \theta u \cdot \nabla \phi)\,dx\,dt + \int_{\R^2}\theta_0(x)\phi(0,x)\,dx =0,
$$
for all $\phi \in C_c^{\infty}([0,T) \times \R^2)$.
\end{definition}
\begin{theorem}\label{LWP H21}Let $\theta_0 \in W^{2,(2,1)}(\R^2)$, then there exists a $T>0$ and a unique weak solution to \eqref{SQG} satisfying
$$
\theta \in L^{\infty}([0,T];W^{2,(2,1)}(\R^2)) \cap C([0,T];W^{1,(2,1)}(\R^2)).
$$
Furthermore, $T \leq \frac{1}{C\norm{\theta_0}{W^{2,(2,1)}}}$ and 
$$
\sup _{0 \leq t \leq T} \|\theta(t)\|_{W^{2,(2,1)}} \leq 2C\|\theta_0\|_{W^{2,(2,1)}}.
$$
\end{theorem}
\begin{proof}
Let $\phi\in C_c^{\infty}(\R^2)$ be nonnegative and satisfy $\norm{\phi}{L^1}=1$. For $\eps>0$, we define $\phi_{\eps}(x)=\eps^{-2}\phi(x/\eps)$. We then consider 
\begin{equation}\tag{$SQG_{\eps}$}
\label{mollified sqg}
\begin{cases}
\partial_t\theta^{\eps} + (u^{\eps}\cdot \nabla)\theta^{\eps}=0 \quad &\text{on } (0,T_{\eps})\times \R^2\\
u^{\eps}=\nabla^{\perp}(-\Delta)^{-\frac{1}{2}}\theta^{\eps} &\text{on } [0,T_{\eps})\times \R^2\\
\theta^{\eps}_{0}=\phi_{\eps} \ast \theta_0 &\text{on } \R^2.
\end{cases}
\end{equation}
Since, for every integer $k\geq 3$, $\theta_0^{\eps}\in H^k(\R^2)$, it follows from \cite{CMT} that there exists $T_{\eps}>0$ and a unique smooth solution $\theta^{\eps}$ to \eqref{mollified sqg} on $[0,T_{\eps}]$ with initial data $\theta^{\eps}_0$. \\
\\
{\bf Uniform bounds on the approximation sequence.} We now show that there exists $T>0$ such that $\theta^\eps$ is uniformly bounded in $L^{\infty}([0,T];W^{2,(2,1)}(\R^2))$. We let $X^{\eps}$ be the flow corresponding to $u^{\eps}$, i.e.
$$
\frac{d}{dt}X^{\eps}(t,x) = u^{\eps}(t,X^{\eps}(t,x)), \quad X^{\eps}(0,x)=x.
$$
Then, for any $i,j \in \{1,2\}$
\begin{equation}
\label{H21 ap eq1}
\begin{split}
&\frac{d}{dt} \theta^{\eps}(t,X^{\eps}(t)) =0,\\
&\frac{d}{dt} (\partial_i\theta^{\eps})(t,X^{\eps}(t))= -(\partial_iu^{\eps}\cdot \nabla \theta^{\eps})(t,X^{\eps}(t)),\\
&\frac{d}{dt} (\partial^2_{i,j}\theta^{\eps})(t,X^{\eps}(t))= -(\partial^{2}_{i,j}u^{\eps}\cdot \nabla \theta^{\eps} + \partial_iu^{\eps}\cdot \partial_j\nabla\theta^{\eps} + \partial_ju^{\eps}\cdot \nabla \partial_i\theta^{\eps})(t,X^{\eps}(t)).
\end{split}
\end{equation}
From Theorem \ref{lorentz embedding two} and Proposition \ref{CZO lorentz}, it follows that
\begin{equation}
\label{gradient bounds for LWP}
\begin{split}
&\|\nabla u^{\eps}\|_{L^{\infty}} + \|\nabla\theta^{\eps}\|_{L^{\infty}} \leq C(\|u^{\eps}\|_{W^{2,(2,1)}} + \|\theta^{\eps}\|_{W^{2,(2,1)}}) \leq C\|\theta^{\eps}\|_{W^{2,(2,1)}}.
\end{split}
\end{equation}
We integrate each equation of \eqref{H21 ap eq1} in time, and we apply the $L^{(2,1)}$ norm to each of the resulting equalities.  We then apply Proposition \ref{Lorentz Holder} and, subsequently,  \eqref{gradient bounds for LWP}, to the right hand side of each inequality. Since Lorentz norms are invariant under measure-preserving transformations, after summing over $i,j$ and applying Proposition \ref{convolution on lorentz}, we conclude that
\begin{equation}
\label{H21 ap eq 2}
\begin{split}
\norm{\theta^{\eps}(t)}{W^{2,(2,1)}} \leq \norm{\theta_0}{W^{2,(2,1)}} + C\int_0^t \norm{\theta^{\eps}(\tau)}{W^{2,(2,1)}}^2\,d\tau.
\end{split}
\end{equation}
In particular, we choose
\begin{equation}
\label{T H21 ap eq 3}
T\leq \frac{1}{2C\norm{\theta_0}{W^{2,(2,1)}}},
\end{equation}
so that, by (\ref{H21 ap eq 2}) and an application of Osgood's lemma (see Chapter $5$ in \cite{Che98}, for example), 
\begin{equation}
\label{H21 uniform bound}
\sup_{\eps > 0}\sup_{0 \leq t \leq T} \norm{\theta^{\eps}(t)}{W^{2,(2,1)}} \leq 2C\norm{\theta_0}{W^{2,(2,1)}},
\end{equation}
hence, by \eqref{gradient bounds for LWP}, the standard continuation criterion ensures that each smooth solution $\theta^{\eps}$ exists throughout $[0,T]$.
Furthermore, from \eqref{H21 uniform bound}, there exists a subsequence $\theta^{\eps_n}$ of $\theta^{\eps}$ which converges weak-$*$ in $L^{\infty}([0,T];W^{2,(2,1)}(\R^2))$.  Let $\theta$ be such that
$$\theta^{\eps_n} \rightharpoonup^* \theta \quad \text{ in }L^{\infty}([0,T];W^{2,(2,1)}(\R^2)).$$ 
\\
{\bf Strong convergence in a lower order norm.} To demonstrate that $\theta$ is a weak solution to ($SQG$), it suffices to show that $\theta^{\eps_n}$ converges to $\theta$ strongly in $L^{\infty}([0,T];L^2(\R^2))$. To this end, we show that $\theta^{\eps_n}$ is Cauchy in $C([0,T];W^{1,(2,1)}(\R^2)).$
For ease of notation set $\overline{\theta}=\theta^{\eps_n}-\theta^{\eps_m}$ and $\overline{u}=u^{\eps_n}-u^{\eps_m}.$ One easily obtains 
$$
\partial_t\overline{\theta} + u^{\eps_n}\cdot \nabla \overline{\theta} = -\overline{u}\cdot \nabla\theta^{\eps_m}.
$$
An argument almost identical to that used to establish \eqref{H21 ap eq 2} yields
\begin{equation*}
\label{H21 Cauchy eq 1}
\begin{split}
\norm{\overline{\theta}(t)}{W^{1,(2,1)}} \leq \norm{\overline{\theta}_0}{W^{1,(2,1)}} + C\int_0^t \norm{\overline{\theta}(\tau)}{W^{1,(2,1)}}&(\norm{\nabla u^{\eps_n}(\tau)}{L^\infty} 
\\ &\ \ +\norm{\nabla \theta^{\eps_m}(\tau)}{L^\infty} + \norm{\theta^{\eps_m}(\tau)}{W^{2,(2,1)}} )d\tau.
\end{split}
\end{equation*}
We apply (\ref{gradient bounds for LWP}), \eqref{H21 uniform bound}, and Gr\"{o}nwall's Lemma to the above estimate, which gives  
\begin{equation}
\label{H21 Cauchy eq 2}
\norm{(\theta^{\eps_n}-\theta^{\eps_m})(t)}{W^{1,(2,1)}} \leq \norm{\theta^{\eps_n}_0 - \theta^{\eps_m}_0}{W^{1,(2,1)}}e^{CT\norm{\theta_0}{W^{2,(2,1)}}},
\end{equation}
showing that $\theta^{\eps_n}$ is Cauchy in $C([0,T];W^{1,(2,1)}(\R^2)).$ Hence, by uniqueness of limits,  
\begin{equation}\label{strongthetaconv}
\theta^{\eps_n} \to \theta \quad\text{strongly in } C([0,T];W^{1,(2,1)}(\R^2)).
\end{equation}
Furthermore, we see that for all $|\alpha| \leq 2$, $D^{\alpha}\theta^{\eps_n} \to D^{\alpha}\theta$ in $\mathcal{D}'([0,T]\times \R^2)$. Hence, by the lower semi-continuity of the Lorentz norms under distributional convergence,  
$$
\norm{\theta}{L^{\infty}([0,T];W^{2,(2,1)})} \leq \liminf_{n \to \infty} \norm{\theta^{\eps_n}}{L^{\infty}([0,T];W^{2,(2,1)})} \leq C.
$$
\\
{\bf The limit $\theta$ is a weak solution.} Finally, we show that $\theta$ is a weak solution to ($SQG$). First note that since for each $n \in \mathbb{N}$, $\theta^{\eps_n}$ is a weak solution, $\text{div }u^{\eps_n}=0$ and for all $\phi \in C_c^{\infty}([0,T)\times \R^2)$, 
\begin{equation*}
    \int_0^T\int_{\R^2} \left(\theta^{\eps_n}\partial_t\phi + \theta^{\eps_n}u^{\eps_n}\cdot \nabla \phi\right)\,dx\,dt + \int_{\R^2}\theta^{\eps_n}_0(x)\phi(0,x)\,dx=0.
\end{equation*}
Since $W^{1,(2,1)}(\R^2)\hookrightarrow L^2(\R^2)$, it follows from (\ref{strongthetaconv}) that 
\begin{equation}
\label{weak soln H21 eq 2}
\left|\int_0^T\int_{\R^2} (\theta^{\eps_n}-\theta)\partial_t\phi \,dx\,dt\right| \leq \sup_{t \leq T} \norm{\theta^{\eps_n}(t)-\theta(t)}{L^2}\int_0^T \norm{\partial_t\phi(t)}{L^2}\,dt \to 0.
\end{equation}
For the nonlinear term, setting $u=R^{\perp}\theta$, we again use (\ref{strongthetaconv}) and boundedness of Riesz operators on $L^2(\R^2)$ to conclude that
\begin{equation*}
\label{weak soln H21 eq3}
\begin{split}
\left| \int_0^T \int_{\R^2} (\theta^{\eps_n}u^{\eps_n}-\theta u)\cdot \nabla \phi\,dx\,dt\right|&\leq \sup_{t \leq T} \norm{\theta^{\eps_n}(t)u^{\eps_n}(t)-\theta(t) u(t)}{L^1}T\norm{\nabla \phi}{L^{\infty}_{t,x}}\\
&\hspace{-4cm} \leq \sup_{t \leq T} \left\{\norm{\theta^{\eps_n}(t)}{L^2} \norm{u^{\eps_n}(t)-u(t)}{L^2}+ \norm{u(t)}{L^2}\norm{\theta^{\eps_n}(t)-\theta(t)}{L^2} \right\}T\norm{\nabla \phi}{L^{\infty}_{t,x}} \to 0.
\end{split}
\end{equation*}
Finally, since $\theta^\eps_0 \to \theta_0$ in $L^2(\R^2)$,
\begin{equation*}
\left|\int_{\R^2} (\theta^{\eps_n}_0-\theta_0)\phi(0,x)\,dx\right| \leq \norm{\theta_0^{\eps_n}-\theta_0}{L^2} \norm{\phi(0)}{L^2} \to 0.
\end{equation*}
Hence $\theta \in L^{\infty}([0,T];W^{2,(2,1)}(\R^2))$ is a weak solution to \eqref{SQG}. \\
\\
{\bf Uniqueness.} Uniqueness follows from an argument identical to that used to show $\theta^{\eps_n}$ is Cauchy in $W^{1,(2,1)}(\R^2)$.
\end{proof}

\end{appendix}
\bibliographystyle{alpha}
\bibliography{sqg}

\end{document}